%% file: MCPDE_arxiv.tex
\documentclass[11pt]{article}

\usepackage{amsmath}
\usepackage{amssymb}
\usepackage{empheq}
\usepackage{txfonts}
\usepackage{graphicx}
\usepackage{caption}
\usepackage{empheq}

\def \E{\mathbb{E}}
\def \R{\mathbb{R}}
\def \M{\mathbb{M}}

\def\Lc{{\cal L}}
\def\x{\times}
\def \Fb{\overline F}
\def\1{{\bf 1}}
\def \I{{\bf I}}
\def \N{\mathbb{N}}
\newcommand{\un}{1\hspace{-1mm}{\rm I}}   % 1

\newtheorem{Theorem}{Theorem}[section]
\newtheorem{proposition}[Theorem]{Proposition}
\newtheorem{remark}[Theorem]{Remark}
\newtheorem{assumption}[Theorem]{Assumption}

\title{Nesting Monte Carlo for high-dimensional Non Linear PDEs}
\author{Xavier Warin \thanks{EDF R\&D \& FiME, Laboratoire de Finance des March\'es de l'Energie, xavier.warin@edf.fr}}

\begin{document}

\maketitle

\begin{abstract}
  A new method based on nesting Monte Carlo is developed
to solve high-dimensional semi-linear PDEs.
Convergence of the method is proved and its convergence rate studied.
Results in high dimension for different kind of non-linearities show its efficiency.
\end{abstract}

\input{sourceTex}
\bibliographystyle{abbrv}
\bibliography{MyBib}

\end{document}

%% file: sourceTex.tex
\section{Introduction}
The resolution of Non Linear PDEs in high dimension is a challenging task due to the so-called "curse of dimensionality".
Deterministic method can't cope with dimensions higher than 4 even using super-computers.
In order to solve problems in higher dimension, effective resolution of Semi-Linear PDEs based on the BSDE approach first proposed by  \cite{pardoux1990adapted} were developed  in \cite{gobet2005regression} and \cite{lemor2006rate}. A lot of literature on the subject has developed  in recent years and the methodology has been extended to solve full non linear PDEs in \cite{fahim2011probabilistic}, \cite{tan2013splitting}.\\
However because the methodology needs some basis functions to project conditional expectation, it faces the curse of dimensionality too by not being able to solve PDEs in dimension higher than six or seven.\\
Recently  two new approaches have emerged  in very high dimensions:
\begin{itemize}
\item The first is based on Deep Learning and uses deep neural networks \cite{han2017overcoming}, \cite{weinan2017deep}, \cite{fujii2017asymptotic}. The method seems to be effective in dimension over 100 but no proof of convergence is currently available and therefore  we don't know its limitations.
\item The second is based on branching methods and is effective for
non linearities polynomial in the solution $u$ and its gradient $Du$. Convergence results are given in \cite{henry2016branching} and numerical results show that the PDEs can be solved in dimension 100. However the authors showed that the variance of the method explodes rapidly when the maturity grows or when the non linearity becomes important: numerical tests confirm that it is in fact  the case.
This methodology has being extended for other non linearities in \cite{bouchard2017numerical} and \cite{bouchard2017numerical2} and the maturity problem is solved but at the price of the introduction of some grids meaning that the "curse of dimensionality" is back.
\item The third is developed in \cite{weinan2017linear}, \cite{weinan2016multilevel}, \cite{hutzenthaler2017multi} with an algorithm based on Picard iterations, multi-level techniques and automatic differentiation permitting to solve some high dimensional PDEs with non linearity in $u$ and $Du$. The convergence of the algorithm is given and  a lot of numerical examples show its efficiency  in high dimension.
\end{itemize}
As for the branching method, the methodology proposed here is based on the Feynman-Kac representation of the PDEs coupled with the randomization of the time step proposed in \cite{henry2016branching}.
This approach is combined with nesting Monte Carlo with a given depth.
Then it is possible to get effective schemes to solve non-linear PDEs.\\ Because a truncation after a given number $m$ of nesting is achieved, the method is biased.\\
For the demonstration of the convergence  of the proposed  schemes it is possible to follow classical approaches as the one used in \cite{rainforth2017opportunities}. It permits to understand how many particles to use at each nesting level and the number of nesting level $m$ to take.\\
Classically the error is composed of a biased term and a variance term. It can be shown than the bias term
goes to zero very quickly with $m$  but to be effective we need to be able to take $m$ below 5 or 6 such that the nesting Monte Carlo can be used. Therefore,
as we will see, a  limitation of the method will be that the maturity  cannot be too large.\\

However the methodology proposed here has a lot of good properties :
\begin{itemize}
\item It is very simple to implement,
\item It needs  a very low memory to run on computers,
\item Its convergence is independent on the dimension $d$ of the problem,
\item It is embarrassingly parallel so it can be run easily on super computers,
\item If the different Lipschitz constants  associated to the non-linearity  are not too large, then the number of particles to take at each level to get a given accuracy is decreasing very fast giving a method very quickly converging.
\end{itemize}

Practically  we will show that the method can be used on a wide set of cases and that we are able for example to solve all test cases proposed in 
\cite{han2017overcoming}, \cite{weinan2017deep} for example.\\
The article has two parts :
\begin{itemize}
\item The first part is devoted to the resolution of the problem  with linearity in $u$. The scheme is given, its convergence studied and numerical results in dimension 6 and 100 show the efficiency of the scheme.
\item The second part is devoted to non linearities in $Du$. Based on automatic differentiation \cite{fournie1999applications}, we give a first scheme and show its convergence. We then introduce a second scheme using ideas in \cite{warin2017variations}.
The second scheme permits to gain little theoretically but numerically we show that is is far more effective than the first proposed. We test the methods on problems with dimensions 10 to 100. 
\end{itemize}

In the sequel we use the classical notation for $Y \in \R^d$ , $||Y||_2 =  \sqrt{\sum_{i=1}^d Y_i^2}$.\\
	Given two matrix $A, B \in \M^d$, denote $A:B := \mbox{Trace}(AB^{\top})$, $\un_d$ is the unit vector of $\R^d$ and $\I_d$ the identity matrix of $\M^d$.\\
 All numerical experiments are achieved on a cluster using 8 nodes with a total of 224 cores and MPI is used for parallelization.  The generation of random numbers  is  achieved using Tina's Random Number Generator Library \cite{bauke2011tina}.
 All computational times are given for a configuration of Intel Xeon CPU E5-2680 v4  2.40GHz (Broadwell).

\section{A first non linear case}
In this section we study the case of a non linearity in $u$ and we aim at solving the PDE  for $ t<T,~x\in\R^d$:
 \begin{flalign}
 \label{eqPDE}
  (-\partial_tu-\Lc u)(t,x)  & = f(t,x,u(t,x)), \nonumber \\
  u_T&=g,
 \end{flalign}
 where
 \begin{flalign}
 \Lc u(t,x) := \mu Du(t,x) +  \frac{1}{2} \sigma \sigma^{\top} \!:\! D^2 u(t,x),
 \label{eq:gen}
 \end{flalign}
so that $\Lc $ is the generator associated to 
 \begin{flalign}
 \hat X_t = x +  \mu t+ \sigma  dW_t,
 \label{eq:sde}
 \end{flalign}
 with $\mu \in \R^d$, $\sigma \in \M^d$ is some constant matrix and $W_t$ a $d$-dimensional Brownian motion.\\
 We will use the following assumptions:
 \begin{assumption}
 \label{ass::lipF}
 $f$ is  uniformly Lipschitz in $u$ with constant $K$ :
 \begin{flalign}
 | f(t, x, y) - f(t, x, w) | \le K |y-w| \quad  \forall t \in \R^{+}, x \in \R^d, (w,y) \in \R^2.
 \end{flalign}
 \end{assumption}

 \begin{assumption}
 \label{ass::uHolder}
 Equation  \eqref{eqPDE} has a solution $u \in C^{1,2}([0,T] \x \R^d)$, such that
 \begin{itemize}
 \item $u$ is $\theta$-H\"older with $\theta \in (0,1]$ in time with constant $\hat K$ :
 \begin{align*}
 | u(t,x) - u(\tilde t, x) | \le \hat K | t- \tilde{t}|^\theta \quad \quad \quad \forall (t,\tilde{t},x) \in [0,T] \times [0,T] \times \R^d,
 \end{align*}
 \item $u(t,x)$ has a quadratic growth in $x$ uniformly in $t$,
 %\item $\E\big( \int_0^T |f(t,X_t, u(t,X_t))| dt\big)< \infty$.
 \end{itemize}
 \end{assumption}
 %\begin{assumption}
 %\label{ass::uBound}
 %$f$ and $g$ satisfies
 %$\sup_{t\in[0,T]} ||f(t,X_t,u(t,X_t))||_2^2   + ||g(X_T)||_2^2 < \infty$
 %\end{assumption}
 In this section 
 $\rho (x)= \lambda e^{- \lambda x}$ is the density of a random variable with exponential  law. \\ Denote $$ \Fb(t):=\int_t^\infty\rho(s)ds = e^{- \lambda t}= 1-F(t),$$
 so that $F$ is the cumulative distribution function of a random variable with density $\rho$.

 \subsection{Idea of the algorithm}
 
 We consider a sequence of i.i.d. random variables $(\tau_{m})_{m \ge 1}$ of density $\rho$. \\
 We consider the sequence defined by: 
 \begin{equation}
 \left\{
 \begin{array}{lll}
 T_0 &  =& 0, \\
 T_{k+1} & =&  (T_{k} + \tau_k) \wedge T.
 \end{array}
 \right.
 \label{eq:T}
 \end{equation}
 We further define $N_T = \inf \{ n | T_{n+1} \ge T \}$.\\
 We also consider a sequence of independent $d$-dimensional Brownian motion $(W^{m}_t)_{m\ge 1}$, which are  independent of $(\tau_{m})_{m \ge 1}$.
 
 Define $W_t =  W^{1}_t$ for all $t \in \big[0, T_{1} \big]$ and
	then for each $k$,
	define 
	\begin{flalign}
    \label{eq:def_Wk}
		W_t ~:=~ W_{T_{k}} + W^{k+1}_{t - T_{k}}, ~~\mbox{for all}~ t \in [T_{k}, T_{k+1}].
	\end{flalign}
	We define an associated diffusion process $(X_t)_{t \in [T_{k}, T_{k+1}]}$ 
	by means of the following SDE:
	\begin{flalign}
    \label{eq:def_Xk}
		X_t ~=~ X_{T_{k}} +  \mu (t-T_k)+ \sigma  W^{k+1}_{t-T_{k}}
		~~t \in [T_{k}, T_{k+1}],~~\P\mbox{-a.s.,}
	\end{flalign}
    with $X_{0} =x $.

 Denoting by $\E_{t,x}$ the expectation operator conditional on  $X_t=x$ at time $t$,  from the Feynman-Kac formula the representation of the solution $u$ valid under assumption \ref{ass::uHolder} is given by:
\begin{flalign*}
 u(0,x) 
 = & 
 \E_{0,x} \Big[ \Fb(T)\frac{g(X_T)}{ \Fb(T)}+\int_0^T \frac{f(t,X_t,u(t,X_t))}{\rho(t)}\rho(t)dt\Big]
 \\ =&
 \E_{0,x} \big[\phi\big(0, T_{1},X_{T_{1} }, u(T_{1} ,X_{T_{1} })\big)\big],
 \end{flalign*}
  \begin{flalign*}
  \phi(s, t,y,z) &:= \frac{\1_{\{t\ge T\}}}{ \Fb(T-s)}g(y)\!+\! \frac{\1_{\{t<T\}}}{\rho(t -s)} f(t,y,z).
 \end{flalign*}
 Recursively we have for $n < N_T$, noting $u_n = u(T_{n}, X_{T_{n}})$ :
 \begin{flalign}
 u_{n} = & \E_{ T_{n}, X_{T_n}}\big[\phi\big(T_n , T_{n+1}, X_{T_{n+1}}, u_{n+1}\big)\big],
 \end{flalign}
 We further consider the truncated operator after $p$ switches :
 \begin{flalign}
 \label{eq:algo}
 u_p^p = & g(X_{T_p}), \nonumber\\
 u_{n}^p = & \E_{ T_{n}, X^n_{T_n}}\big[\phi\big(T_n, T_{n+1}, X_{T_{n+1}}, u_{n+1}^p\big)\big], \quad n < p, \quad  \mbox{  defined if  } T_{n} < T
 \end{flalign}
 The goal of this section is to study the underlying algorithm when
 the resolution of  equation \eqref{eq:algo} is achieved by nesting Monte Carlo. Starting from the ideas used in \cite{henry2016branching} we propose
 a nesting algorithm calculating all $u_{n}^p$ by Monte Carlo.
 We have to show the bias associated to the algorithm goes to  zero and that the global variance induced is controlled.
 In order to get a useful algorithm, we have to show that the bias goes to zero very quickly so that the number of switches to take is low: indeed it is well known that nesting Monte Carlo is subject to an explosion of the computer time. We will show that for many useful cases it is an effective approach.
 
 \subsection{Estimator and global error}
 \label{sec::estim}
Let set $p \in \N^{+}$.
For $(N_0,.., N_{p-1}) \in \N^p$, we introduce  the sets of i-tuple
$Q_i = \{ k=(k_1, ...,k_i)\}$ for $i \in \{1,..,p\}$ where all components $k_j \in [1, N_{j-1}]$. 
Besides we define $Q^p= \cup_{i=1}^p Q_i$.\\
For $k=(k_1, ...,k_i) \in Q_i$ we introduce the set $\tilde Q(k) = \{ l =(k_1,..,k_i, m )/  m \in \{1,.., N_i\} \} \subset Q_{i+1}$. By convention $\tilde Q(\emptyset) = \{ l =( m )/  m \in \{1,.., N_0\} \} =  Q_{1}.$\\
We define the following sequence $\tau_{k}$  of switching increments always i.i.d. random variables with density $\rho$ 
 for $k\in Q^p$, and  
 a sequence of independent $d$-dimensional Brownian motion $(\bar W^{k})$, which are  independent of the $(\tau_{k})_{k \in Q^p}$.
 Let us define the switching dates:
 \begin{equation}
 \left\{
 \begin{array}{lll}
 T_{(j)} &  =&  \tau_{(j)} \wedge T,  j \in \{ 1,., N_{0}\} \\
 T_{ \tilde k} & =&  ( T_{k} + \tau_{ \tilde k}) \wedge T,  k =( k_1,..k_i) \in Q_i, \tilde k \in \tilde Q(k)
 \end{array}
 \right.
 \end{equation}
 
We define an associated diffusion process $(X^{\tilde k}_t)_{t\geq 0}$
by means of the following SDE
\begin{flalign}
  X^{(i)}_t ~=~&  X_{0}^\emptyset + \mu t  
		+\sigma  \bar W^{(i)}_t,
		~~t \in [0, T_{(i)}], i=1, N_0 \nonumber \\
  X^{\tilde k}_t ~=~& X^{k}_{T_{k}} + \mu  (t-T_k)
		+ \sigma  \bar W^{\tilde k}_{t-T_k}, \mbox{ for } \tilde{k} \in \tilde Q(k), 
		~~t \in [T_{k}, T_{\tilde k}],~~\P\mbox{-a.s.,} 
        \label{eq:diff}
\end{flalign}
with $X_{0}^\emptyset =x $.

We consider the estimator defined by:
\begin{align}
\label{eq:estim1}
\left\{
 \begin{array}{ll}
\bar u_\emptyset^p = & \frac{1}{N_0} \sum_{j=1}^{N_0} \phi\big( 0,T_{(j)}, X^{(j)}_{T_{(j)}}, \bar u_{(j)}^p\big),  \\
\bar u_{k}^p = &  \frac{1}{N_i} \sum_{\tilde k \in \tilde Q(k)} \phi\big( T_k, T_{\tilde k},  X^{\tilde k}_{T_{\tilde k}}, \bar u_{\tilde k}^p\big), \\
& \quad \mbox{ for } k =( k_1,..k_i) \in Q_i, 1 < i <p , T_{k} < T, \\
\bar u_{\tilde k}^p = & g(X^{\tilde k}_{T_{\tilde k}})  \quad \mbox{ for } \tilde k \in Q_p.
\end{array}
\right.
\end{align}

Note that in the case where $T_{\tilde{k}} = T$ then $\phi\big( T_k, T_{\tilde k},  X^{\tilde k}_{T_{\tilde k}}, \bar u_{\tilde k}^p\big)$ is independent of $\bar u_{\tilde k}^p$ so that the recursion is stopped.

\begin{proposition}
\label{theo1}
Under assumption \ref{ass::lipF}, \ref{ass::uHolder} , we have the following error given by the estimator \eqref{eq:estim1} :
\begin{align*}
\E \big((\bar u_{\emptyset}^p - u(0,x))^2 \big) \le   \prod_{i=1}^p  (1+\frac{8}{N_{i-1}}) \frac{K^{2p} e^{\lambda T} }{\lambda^{p}} T^{2\theta} \hat{K}^2 \frac{T^p}{p\Gamma(p)}+ \nonumber \\
\sum_{i=0}^{p-1} \frac{K^{2i}}{N_{i}} \prod_{j=1}^i  (1+\frac{8}{N_{j-1}})
\frac{ T^{i} e^{\lambda T}}{\lambda^{i}}
\kappa_i
\end{align*}
with
\begin{align*}
\kappa_{i} = \big( \frac{4T }{\lambda (i+1)!} \sup_{t\in[0,T]} \E\big( f(t,X_t,u(t,X_t))^2 \big)	+    2 \E\big( g(X_T)^2 \big)  (\frac{ 1_{i>0}}{i \Gamma(i)} + 1_{i=0})  \big)
\end{align*}
and  $\Gamma(s)= \int_0^\infty t^{s-1} e^{-t} dt$ is the gamma function.
\end{proposition}
\begin{proof}
First notice that due to the Lipschitz property in assumption \ref{ass::lipF} and the growth assumption on $u$, $\E\big( \int_0^T |f(t,X_t, u(t,X_t))| dt\big)< \infty$.
Then notice that under assumption \ref{ass::uHolder}, the solution $u$ of \eqref{eqPDE} satisfies  a Feynman-Kac relation  (see an adaptation of proposition 1.7 in \cite{touzi2012optimal} ) so that
for all $k \in Q_i$,  and $\forall \tilde k \in \tilde{Q}(k)$,
\begin{align}
u(T_k, X^k_{T_k}) =& \E_{ T_{k}, X^k_{T_k}} \big[\phi\big(T_k , T_{\tilde k}, X_{T_{\tilde k}^{\tilde{k}}}, u(T_{\tilde k},X_{T_{\tilde k}^{\tilde{k}}} ))\big],
\label{eqq:realU}
\end{align}
Then for $k \in Q_i$, $i <p$, let us  define:
\begin{flalign}
\label{eq:err}
E_k := & \E_{T_{k},X^{k}_{T_{k}} }\big((\bar u_{k}^p - u(T_k,X^k_{T_k}))^2 1_{T_k <T} \big) \nonumber\\
= &  V_k  + B_k^2
\end{flalign}
where we note $B_k$ the bias error for index $k$ as:
\begin{flalign}
B_k :=  \big( \E_{T_{k},X^{k}_{T_{k}}} (\bar u_{k}^p) -u(T_k,X^k_{T_k}) \big) 1_{T_k <T} ,
\end{flalign}
and the variance term $v_k$ of the estimator:
\begin{flalign}
V_k := \E_{T_{k} ,X^{k}_{T_{k}} }\big( 1_{T_k <T}(\bar u_{k}^p- \E_{T_{k} ,X^{k}_{T_{k}} }(\bar u_{k}^p))^2).
\label{eq::V}
\end{flalign}
Let us begin with the variance term.\\
Note that using the equation \eqref{eqq:realU} and the $\bar u_{k}^p$ definition given by equation \eqref{eq:estim1}:
\begin{flalign*}
%\label{eq:p1}
(\bar u_{k}^p - \E_{T_k, X^{k}_{T_{k}}}(\bar u_{k}^p) )^2 1_{T_k <T}  & =  \big[\frac{1}{N_i} \sum_{\tilde k \in \tilde Q(k)}  ( 1_{T_{\tilde k < T}} \frac{f(T_{\tilde k}, X^{\tilde k}_{T_{\tilde k}}, \bar u_{\tilde k}^p)}{\rho(\tau_{\tilde k})} -\E_{T_k,X^{k}_{T_{k}}}(1_{T_{\tilde k < T}} \frac{f(T_{\tilde k}, X^{\tilde k}_{T_{\tilde k}},\bar u_{\tilde k}^p)}{\rho(\tau_{\tilde k})})) + \nonumber \\
&  \frac{1}{N_i}   \sum_{\tilde k \in \tilde Q(k)}  ( 1_{T_{\tilde k\ge T}} \frac{g(X^{\tilde k}_T)}{ \Fb(T-T_k)} - \E_{T_k,X^{k}_{T_{k}}}(1_{T_{\tilde k\ge T}} \frac{g(X^{\tilde k}_T)}{ \Fb(T-T_k)}))\big]^2
\end{flalign*}
so that
\begin{flalign*}
  V_k \le &2 \E_{T_{k} ,X^{k}_{T_{k}} }\big[\big( \frac{1}{N_i} \sum_{\tilde k \in \tilde Q(k)}  (1_{T_{\tilde k < T}} \frac{f(T_{\tilde k}, X^{\tilde k}_{T_{\tilde k}}, \bar u_{\tilde k}^p)}{\rho(\tau_{\tilde k})} -\E_{T_k,X^{k}_{T_{k}}}(1_{T_{\tilde k < T}} \frac{f(T_{\tilde k}, X^{\tilde k}_{T_{\tilde k}},\bar u_{\tilde k}^p)}{\rho(\tau_{\tilde k})}))\big)^2\big] + \nonumber \\
& 2 \E_{T_{k} ,X^{k}_{T_{k}} }\big[\big(\frac{1}{N_i}   \sum_{\tilde k \in \tilde Q(k)}  ( 1_{T_{\tilde k\ge T}} \frac{g(X^{\tilde k}_T)}{ \Fb(T-T_k)} - \E_{T_k,X^{k}_{T_{k}}}(1_{T_{\tilde k\ge T}} \frac{g(X^{\tilde k}_T)}{ \Fb(T-T_k)}))\big)^2\big]
\end{flalign*}
and using that for independent random variables $x_i$ 
\begin{flalign*}
\E[(\frac{1}{N} \sum_{i=1}^N (x_i-E(x_i)))^2]= \frac{1}{N^2} \sum_{i=1}^N \E[(x_i-\E(x_i))^2]
\end{flalign*}
\begin{flalign}
  \label{eq:p1}
  V_k \le  &\frac{2}{N_i^2} \sum_{\tilde k \in \tilde Q(k)} \E_{T_{k} ,X^{k}_{T_{k}} }\big[ \big(  1_{T_{\tilde k < T}} \frac{f(T_{\tilde k}, X^{\tilde k}_{T_{\tilde k}}, \bar u_{\tilde k}^p)}{\rho(\tau_{\tilde k})} -\E_{T_k,X^{k}_{T_{k}}}(1_{T_{\tilde k < T}} \frac{f(T_{\tilde k}, X^{\tilde k}_{T_{\tilde k}},\bar u_{\tilde k}^p)}{\rho(\tau_{\tilde k})})\big)^2 \big] + \nonumber \\
& \frac{2}{N_i^2}  \sum_{\tilde k \in \tilde Q(k)}   \E_{T_{k} ,X^{k}_{T_{k}} }\big[\big(  1_{T_{\tilde k\ge T}} \frac{g(X^{\tilde k}_T)}{ \Fb(T-T_k)} - \E_{T_k,X^{k}_{T_{k}}}(1_{T_{\tilde k\ge T}} \frac{g(X^{\tilde k}_T)}{ \Fb(T-T_k)})\big)^2 \big].
\end{flalign}
Developing
\begin{flalign*}
A = & 1_{T_{\tilde k < T}} \frac{f(T_{\tilde k}, X^{\tilde k}_{T_{\tilde k}}, \bar u_{\tilde k}^p)}{\rho(\tau_{\tilde k})} -\E_{T_k,X^k_{T_k}}(1_{T_{\tilde k < T}} \frac{f(T_{\tilde k}, X^{\tilde k}_{T_{\tilde k}},\bar u_{\tilde k}^p)}{\rho(\tau_{\tilde k})})
\\& = (1_{T_{\tilde k < T}} \frac{f(T_{\tilde k}, X^{\tilde k}_{T_{\tilde k}}, \bar u_{\tilde k}^p) -f(T_{\tilde k}, X^{\tilde k}_{T_{\tilde k}}, u(T_{\tilde k},X^{\tilde k}_{T_{\tilde k}})) }{\rho(\tau_{\tilde k})}) + \\
& ( 1_{T_{\tilde k < T}} \frac{f(T_{\tilde k}, X^{\tilde k}_{T_{\tilde k}},u(T_{\tilde k},X^{\tilde k}_{T_{\tilde k}}) )}{\rho(\tau_{\tilde k})} - \E_{T_k,X^k_{T_k}}(1_{T_{\tilde k < T}} \frac{f(T_{\tilde k}, X^{\tilde k}_{T_{\tilde k}}, \bar u_{\tilde k}^p)}{\rho(\tau_{\tilde k})}))
\end{flalign*}
that we inject in \eqref{eq:p1} so  that using the relation $\E(|x-\E(x)|^2) \le \E(x^2)$,  $V_k$  is bounded following:
\begin{flalign}
\label{eq:var}
V_k \le  & 4   \frac{1}{N_i^2} \sum_{\tilde k \in \tilde Q(k)} \E_{T_k,X^k_{T_k} } \big( 1_{T_{\tilde k < T}} (\frac{f(T_{\tilde k}, X^{\tilde k}_{T_{\tilde k}}, \bar u_{\tilde k}^p) -f(T_{\tilde k}, X^{\tilde k}_{T_{\tilde k}},u(T_{\tilde k},X^{\tilde k}_{T_{\tilde k}}))}{\rho(\tau_{\tilde k})} )^2     \big) + \nonumber\\
& 4 \frac{1}{N_i^2} \sum_{\tilde k \in \tilde Q(k)} \E_{T_k,X^k_{T_k} }\big( ( 1_{T_{\tilde k < T}} \frac{f(T_{\tilde k}, X^{\tilde k}_{T_{\tilde k}}, u(T_{\tilde k},X^{\tilde k}_{T_{\tilde k}}))}{\rho(\tau_{\tilde k}))}-\E_{T_k,X^k_{T_k}}(1_{T_{\tilde k < T}} \frac{f(T_{\tilde k}, X^{\tilde k}_{T_{\tilde k}}, \bar u_{\tilde k}^p)}{\rho(\tau_{\tilde k})}) )^2 \big) + \nonumber \\
& 2 \frac{1}{N_i^2} \sum_{\tilde k \in \tilde Q(k)} \E_{T_k,X^k_{T_k} }\big( 1_{T_{\tilde k\ge T}} \frac{g(X^{\tilde k}_T)^2}{ \Fb(T-T_k)^2} \big).
\end{flalign}
Using that for $X$, $Y$ random variables $\E\big( (X-E(Y))^2\big) = \E\big( (X-E(X))^2\big) + \big(\E(X-Y)\big)^2$, 
then  for $\tilde k  \in \tilde Q(k)$:
\begin{flalign}
I=&\E_{T_k,X^k_{T_k} }\big( ( 1_{T_{\tilde k < T}} \frac{f(T_{\tilde k}, X^{\tilde k}_{T_{\tilde k}}, u(T_{\tilde k},X^{\tilde k}_{T_{\tilde k}}))}{\rho(\tau_{\tilde k})} -\E_{T_k,X^k_{T_k}}(1_{T_{\tilde k < T}} \frac{f(T_{\tilde k}, X^{\tilde k}_{T_{\tilde k}}, \bar u_{\tilde k}^p)}{\rho(\tau_{\tilde k})}))^2 \big) \nonumber \\
= &\E_{T_k,X^k_{T_k} }\big( ( 1_{T_{\tilde k < T}} \frac{f(T_{\tilde k}, X^{\tilde k}_{T_{\tilde k}}, u(T_{\tilde k},X^{\tilde k}_{T_{\tilde k}}))}{\rho(\tau_{\tilde k})} -\E_{T_k,X^k_{T_k}}(1_{T_{\tilde k < T}} \frac{f(T_{\tilde k}, X^{\tilde k}_{T_{\tilde k}}, u(T_{\tilde k},X^{\tilde k}_{T_{\tilde k}}))}{\rho(\tau_{\tilde k})}))^2 \big) + \nonumber\\
&(\E_{T_k,X^k_{T_k}}(1_{T_{\tilde k < T}} \frac{f(T_{\tilde k}, X^{\tilde k}_{T_{\tilde k}}, u(T_{\tilde k},X^{\tilde k}_{T_{\tilde k}}))}{\rho(\tau_{\tilde k})}) - \E_{T_k,X^k_{T_k}}(1_{T_{\tilde k < T}} \frac{f(T_{\tilde k}, X^{\tilde k}_{T_{\tilde k}}, \bar u_{\tilde k}^p))}{\rho(\tau_{\tilde k})}))^2 \nonumber\\
\le  &\E_{T_k,X^k_{T_k} }\big( ( 1_{T_{\tilde k < T}} \frac{f(T_{\tilde k}, X^{\tilde (k}_{T_{\tilde k}}, u(T_{\tilde k},X^{\tilde k}_{T_{\tilde k}}))}{\rho(\tau_{\tilde k})} -\E_{T_k,X^k_{T_k}}(1_{T_{\tilde k < T}} \frac{f(T_{\tilde k}, X^{\tilde k}_{T_{\tilde k}}, u(T_{\tilde k},X^{\tilde k}_{T_{\tilde k}}))}{\rho(\tau_{\tilde k})}))^2 \big) + \nonumber\\
&\E_{T_k,X^k_{T_k} } \big( 1_{T_{\tilde k < T}} (\frac{f(T_{\tilde k}, X^{\tilde k}_{T_{\tilde k}}, \bar u_{\tilde k}^p) -f(T_{\tilde k}, X^{\tilde k}_{T_{\tilde k}},u(T_{\tilde k},X^{\tilde k}_{T_{\tilde k}}))}{\rho(\tau_{\tilde k})} )^2     \big) 
\label{eq:I}
\end{flalign}
where the last inequality is obtained by Jensen.\\
Plugging \eqref{eq:I}  in \eqref{eq:var} and using  the relation
$\E(|x-\E(x)|^2) \le \E(x^2)$  we get  :
\begin{flalign}
V_k \le  & 8   \frac{1}{N_i^2} \sum_{\tilde k \in \tilde Q(k)} \E_{T_k,X^k_{T_k} } \big( 1_{T_{\tilde k < T}} (\frac{f(T_{\tilde k}, X^{\tilde k}_{T_{\tilde k}}, \bar u_{\tilde k}^p) -f(T_{\tilde k}, X^{\tilde k}_{T_{\tilde k}},u(T_{\tilde k},X^{\tilde k}_{T_{\tilde k}}))}{\rho(\tau_{\tilde k})} )^2     \big) + \nonumber\\
& 4 \frac{1}{N_i^2} \sum_{\tilde k \in \tilde Q(k)} \E_{T_k,X^k_{T_k} }\big(  1_{T_{\tilde k < T}} (\frac{f(T_{\tilde k}, X^{\tilde k}_{T_{\tilde k}}, u(T_{\tilde k},X^{\tilde k}_{T_{\tilde k}}))}{\rho(\tau_{\tilde k})} )^2\big) +  \nonumber \\
&2 \frac{1}{N_i^2} \sum_{\tilde k \in \tilde Q(k)} \E_{T_k,X^k_{T_k} }\big( 1_{T_{\tilde k\ge T}} \frac{g(X^{\tilde k}_T)^2}{ \Fb(T-T_k)^2} \big).
\label{eq:V2}
\end{flalign}

Using the Lipschitz property of $f$ and the tower property :
\begin{flalign}
V_k \le  & 8   \frac{1}{N_i^2} \sum_{\tilde k \in \tilde Q(k)} \E_{T_k,X^k_{T_k} } \big(  (\frac{K}{\rho(\tau_{\tilde k})} )^2  E_{\tilde{k}}   \big) + \nonumber\\
& 4 \frac{1}{N_i^2} \sum_{\tilde k \in \tilde Q(k)} \E_{T_k,X^k_{T_k} }\big(  1_{T_{\tilde k < T}} (\frac{f(T_{\tilde k}, X^{\tilde k}_{T_{\tilde k}}, u(T_{\tilde k},X^{\tilde k}_{T_{\tilde k}}))}{\rho(\tau_{\tilde k})} )^2\big) +  \nonumber \\
&2 \frac{1}{N_i^2} \sum_{\tilde k \in \tilde Q(k)} \E_{T_k,X^k_{T_k} }\big( 1_{T_{\tilde k\ge T}} \frac{g(X^{\tilde k}_T)^2}{ \Fb(T-T_k)^2} \big).
\label{eq:V3}
\end{flalign}
Now we take care of the bias term. 
\begin{flalign}
B_k^2 =& (\frac{1}{N_i} \sum_{\tilde k \in \tilde Q(k)}   E_{T_k,X^{k}_{T_{k}}} \big(
1_{T_{\tilde k < T}}   \big( \frac{f(T_{\tilde{k}},X^{ \tilde k}_{T_{\tilde k}},\bar u_{\tilde k}^p)}{\rho(\tau_{\tilde k})} - \frac{f(T_{\tilde{k}},X^{\tilde k}_{T_{\tilde k}},u(T_{\tilde{k}},X^{ \tilde k}_{T_{\tilde k}}))}{\rho(\tau_{\tilde k})}   \big) \big) \big)^2
\label{eq::bk2}
\end{flalign}
Using the fact that all expectations  for the $\tilde k$ are the same, we get :
\begin{flalign*}
B_k^2 =& \frac{1}{N_i} \sum_{\tilde k \in \tilde Q(k)}  E_{T_k,X^{k}_{T_{k}}} \big(
1_{T_{\tilde k < T}}   \big( \frac{f(T_{\tilde{k}},X^{ \tilde k}_{T_{\tilde k}},\bar u_{\tilde k}^p)}{\rho(\tau_{\tilde k})} - \frac{f(T_{\tilde{k}},X^{\tilde k}_{T_{\tilde k}},u(T_{\tilde{k}},X^{ \tilde k}_{T_{\tilde k}}))}{\rho(\tau_{\tilde k})}   \big)\big)^2.
%\label{eq::bk2}
\end{flalign*}
Then using Jensen :
\begin{flalign*}
B_k^2 \le & \frac{1}{N_i} \sum_{\tilde k \in \tilde Q(k)}   E_{T_k,X^{k}_{T_{k}}}\big(
1_{T_{\tilde k < T}}   \big( \frac{f(T_{\tilde{k}},X^{ \tilde k}_{T_{\tilde k}},\bar u_{\tilde k}^p)}{\rho(\tau_{\tilde k})} - \frac{f(T_{\tilde{k}},X^{\tilde k}_{T_{\tilde k}},u(T_{\tilde{k}},X^{ \tilde k}_{T_{\tilde k}}))}{\rho(\tau_{\tilde k})}   \big)^2 \big).
\end{flalign*}
So that using the Lipschitz property of $f$ and the tower property :
\begin{flalign}
\label{eq:bk}
B_k^2 \le & \frac{1}{N_i} \sum_{\tilde k \in \tilde Q(k)} E_{T_k,X^{k}_{T_{k}}} \big(  \frac{K^2}{\rho(\tau_{\tilde k})^2} E_{\tilde k} \big).
\end{flalign}
Plugging  \eqref{eq:V3} and \eqref{eq:bk} in \eqref{eq:err} we get:
\begin{align}
E_k \le& \frac{1}{N_i} (1+\frac{8}{N_i}) \sum_{\tilde k \in \tilde Q(k)} E_{T_k,X^{k}_{T_{k}}} \big(  \frac{K^2}{\rho(\tau_{\tilde k})^2} E_{\tilde k} \big) + \nonumber
\\
& 4 \frac{1}{N_i^2} \sum_{\tilde k \in \tilde Q(k)} \E_{T_k,X^k_{T_k} }\big(  1_{T_{\tilde k < T}} (\frac{f(T_{\tilde k}, X^{\tilde{k}}_{T_{\tilde k}}, u(T_{\tilde k},X^{\tilde k}_{T_{\tilde k}}))}{\rho(\tau_{\tilde k})} )^2\big) +  \nonumber \\
&2 \frac{1}{N_i^2} \sum_{\tilde k \in \tilde Q(k)} \E_{T_k,X^k_{T_k} }\big( 1_{T_{\tilde k\ge T}} \frac{g(X^{\tilde k}_T)^2}{ \Fb(T-T_k)^2} \big).
\label{eq:ErrLess}
\end{align}
We can iterate to get $E_\emptyset$ using the tower property
\begin{align}
\label{eq:E0a}
E_\emptyset \le & \prod_{i=1}^p \frac{1}{N_{i-1}} (1+\frac{8}{N_{i-1}}) \sum_{\tilde k^1 \in \tilde Q(\emptyset)} ... \sum_{\tilde k^p \in \tilde Q(\tilde{k}^{p-1})}  \E[ \frac{K^{2p}}{\prod_{j=1}^p \rho(\tau_{\tilde{k}^j})^2} E_{\tilde{k}^p}] + \nonumber \\
&  \sum_{i=0}^{p-1} \frac{K^{2i}}{N_{i}^2} \prod_{j=1}^i \frac{1}{N_{j-1}} (1+\frac{8}{N_{j-1}})
\sum_{\tilde k^1 \in \tilde Q(\emptyset)} ... \sum_{\tilde k^{i+1} \in \tilde Q(\tilde{k}^{i})}  \E \big [ 1_{T_{\tilde{k}^{i+1} < T}} 
\frac{4 f(T_{\tilde k^{i+1}}, X^{\tilde k ^{i+1}}_{T_{\tilde k^{i+1}}}, u(T_{\tilde k^{i+1}},X^{\tilde k^{i+1}}_{T_{\tilde k^{i+1}}}))^2}{ \prod_{j=1}^{i+1} \rho(\tau_{\tilde k^j})^2}  +\nonumber \\
& 1_{T_{\tilde{k}^{i+1} > T}} 1_{T_{\tilde{k}^{i} < T}} \frac{2 g(X^{\tilde k ^{i+1}}_{T})^2}{\prod_{j=1}^{i} \rho(\tau_{\tilde k^j})^2 \bar F(T-T_{\tilde{k}^{i}})^2} \big ]
\end{align}
where $ E_{\tilde{k}^p} = 1_{T_{\tilde k^p < T}}  (g(X_{T_{\tilde{k}^p}}^{\tilde{k}^p})- u(T_{\tilde{k}^p},X_{T_{\tilde{k}^p}}^{\tilde{k}^p}))^2$.\\
We now bound the two terms in the last summation. 
Using the fact that $\rho$ corresponds to the density of an exponential law:
\begin{align}
D := & \E \big [ 1_{T_{\tilde{k}^{i+1} < T}} 
\frac{f(T_{\tilde k^{i+1}}, X^{\tilde k ^{i+1}}_{T_{\tilde k^{i+1}}}, u(T_{\tilde k^{i+1}},X^{\tilde k^{i+1}}_{T_{\tilde k^{i+1}}}))^2}{ \prod_{j=1}^{i+1} \rho(\tau_{\tilde k^j})^2} \big] \nonumber \\
\le & \E \big [ \frac{1_{T_{\tilde{k}^{i+1} < T}}}{\prod_{j=1}^{i+1} \rho(\tau_{\tilde k^j})} \frac{f(T_{\tilde k^{i+1}}, X^{\tilde k ^{i+1}}_{T_{\tilde k^{i+1}}}, u(T_{\tilde k^{i+1}},X^{\tilde k^{i+1}}_{T_{\tilde k^{i+1}}}))^2}{ \prod_{j=1}^{i+1} \rho(\tau_{\tilde k^j})} \big] \nonumber \\
\le & \frac{e^{\lambda T}}{\lambda^{i+1}} \E \big [ 1_{T_{\tilde{k}^{i+1} < T}} \frac{f(T_{\tilde k^{i+1}}, X^{\tilde k ^{i+1}}_{T_{\tilde k^{i+1}}}, u(T_{\tilde k^{i+1}},X^{\tilde k^{i+1}}_{T_{\tilde k^{i+1}}}))^2}{ \prod_{j=1}^{i+1} \rho(\tau_{\tilde k^j})} \big] \nonumber \\
= & \frac{e^{\lambda T}}{\lambda^{i+1}}  \int_0^T \int_{0}^{T-t^1}..\int_0^{T-t^1-..-t^i} \E[f( \sum_{j=1}^{i+1} t_j, X_{\sum_{j=1}^{i+1} t^j}^{\tilde k ^{i+1}},
u(\sum_{j=1}^{i+1} t^j, X_{\sum_{j=1}^{i+1} t^j}^{\tilde k ^{i+1}}))^2 ] dt^1.. dt^{i+1} \nonumber \\
\le & \frac{e^{\lambda T}}{\lambda^{i+1}}  \frac{T^{i+1}}{(i+1)!} \sup_{t\in[0,T]} \E\big( f(t,X_t,u(t,X_t))^2),	
\label{eq:D}
\end{align}
where $\E\big( f(t,X_t,u(t,X_t))^2) <	\infty$ due to the Lipschitz condition on $f$ and the quadratic  growth of $u$.\\ 
We know deal with the last term using the fact that distribution of $X_T$ is independent of the switching dates :
\begin{align}
H_i:=& \E \big [1_{T_{\tilde{k}^{i+1} > T}} 1_{T_{\tilde{k}^{i} < T}} \frac{g(X^{\tilde k ^{i+1}}_{T_{\tilde k^{i+1}}})^2}{\prod_{j=1}^{i} \rho(\tau_{\tilde k^j})^2 \bar F(T-T_{\tilde{k}^{i}})^2} \big ] \nonumber \\
=  & \E \big [1_{T_{\tilde{k}^{i+1} > T}} \frac{1_{T_{\tilde{k}^{i}} < T}}{\prod_{j=1}^{i} \rho(\tau_{\tilde k^j})^2 \bar F(T-T_{\tilde{k}^{i}})^2} \big]  \E \big( g(X_{T})^2 \big)  \nonumber \\
= &   \E \big [1_{T_{\tilde{k}^{i+1} > T}} 1_{T_{\tilde{k}^{i}} < T} \frac{e^{2\lambda T_{\tilde k^i}}}{\lambda^{2i}} e^{2\lambda (T-T_{\tilde k^i})} \big]  \E\big( g(X_T)^2 \big)  \nonumber\\
= &  \frac{e^{2\lambda T}}{\lambda^{2i}} \E \big [1_{T_{\tilde{k}^{i+1} > T}} 1_{T_{\tilde{k}^{i}} < T} \big] \E\big( g(X_T)^2 \big)   \nonumber \\
= & \frac{e^{2\lambda T}}{\lambda^{2i}} \E\big( g(X_T)^2 \big)   \int_0^T \lambda^{i} x^{i-1} \frac{e^{-\lambda x}} {\Gamma(i)} \int_{T-x}^\infty \lambda e^{-\lambda y} dy dx \nonumber\\
=&  \frac{e^{\lambda T}}{\lambda^{i}} \frac{ T^i}{i \Gamma(i)} \E\big( g(X_T)^2 \big)  \quad \mbox{ for } i > 0, \nonumber \\
H_0 = & \E \big [1_{T_{\tilde{k}^{1}} > T}  \frac{g(X^{\tilde k ^{i+1}}_{T_{\tilde k^{i+1}}})^2}{\bar F(T)^2} \big ] \nonumber \\
    = & e^{\lambda T} 	\E\big( g(X_T)^2 \big).
\label{eq:H}
\end{align}
where we have used the fact that $T_{\tilde{k}^{i}}$ follows a gamma law with density
$\lambda^{i} x^{i-1} \frac{e^{-\lambda x}}{\Gamma(i)}$.\\
At last using the H\"older property of $u$ with respect to $t$:
\begin{align}
\label{eq:F}
F:= & \E[ 1_{T_{\tilde{k}^p < T}} \frac{K^{2p}}{\prod_{j=1}^p \rho(\tau_{\tilde{k}^j})^2} E_{\tilde{k}^p}] \nonumber \\
=  &\E[ 1_{T_{\tilde{k}^p < T}} \frac{K^{2p}}{\prod_{j=1}^p \rho(\tau_{\tilde{k}^j})^2} ( g(X^{\tilde{k}^p}_{T_{\tilde{ k}^p}})-u(T_{\tilde{ k^p}},X^{\tilde{k}^p}_{T_{\tilde{ k}^p}}))^2 ] \nonumber \\
\le & K^{2p} T^{2\theta} \hat{K}^2 \E[ 1_{T_{\tilde{k}^p < T}} \frac{1}{\prod_{j=1}^p \rho(\tau_{\tilde{k}^j})^2}] \nonumber \\
= & \frac{K^{2p}}{\lambda^{2p}} T^{2\theta} \hat{K}^2 \E[ 1_{T_{\tilde{k}^p} < T} e^{ 2\lambda T_{\tilde{k}^p }}] = \frac{K^{2p}}{\lambda^{2p}} T^{2\theta} \hat{K}^2 \int_0^T e^{\lambda x} \frac{\lambda^p x^{p-1}}{\Gamma(p)} dx  \nonumber\\
\le & \frac{K^{2p} e^{\lambda T}}{\lambda^{p}} T^{2\theta} \hat{K}^2 \frac{T^p}{p \Gamma(p) }.
\end{align}
Plugging \eqref{eq:D},\eqref{eq:H} and \eqref{eq:F} in \eqref{eq:E0a} 
\begin{align}
E_\emptyset \le & \prod_{i=1}^p \frac{1}{N_{i-1}} (1+\frac{8}{N_{i-1}}) \sum_{\tilde k^1 \in \tilde Q(\emptyset)} ... \sum_{\tilde k^p \in \tilde Q(\tilde{k}^{p-1})} \frac{K^{2p} e^{\lambda T}}{\lambda^{p}} T^{2\theta} \hat{K}^2 \frac{T^p}{p\Gamma(p)} + \nonumber \\
&  \sum_{i=0}^{p-1} \frac{K^{2i}}{N_{i}} \prod_{j=1}^i  (1+\frac{8}{N_{j-1}})
\frac{T^{i} e^{\lambda T}}{\lambda^{i}}
\big( \frac{4 T  }{\lambda (i+1)!}  \sup_{t\in[0,T]} \E\big( f(t,X_t,u(t,X_t))^2)	 +   \nonumber\\
& 2 \E\big( g(X_T)^2 \big)   (\frac{ 1_{i>0}}{i \Gamma(i)} + 1_{i=0})   \big) \nonumber \\
= &  \prod_{i=1}^p  (1+\frac{8}{N_{i-1}}) \frac{K^{2p} e^{\lambda T}}{\lambda^{p}} T^{2\theta} \hat{K}^2 \frac{T^p}{p\Gamma(p)}+ \nonumber \\
&  \sum_{i=0}^{p-1} \frac{K^{2i}}{N_{i}} \prod_{j=1}^i  (1+\frac{8}{N_{j-1}})
\frac{ T^{i} e^{\lambda T}}{\lambda^{i}}
\big( \frac{4 T}{\lambda (i+1)!}  \sup_{t\in[0,T]} \E\big( f(t,X_t,u(t,X_t))^2)	 + \nonumber \\
& 2 \E\big( g(X_T)^2 \big)  (\frac{ 1_{i>0}}{i \Gamma(i)} + 1_{i=0}) \big)
\end{align}
which is the desired result.
\qquad \end{proof}
Observe that it is necessary to solve accurately  each inner iteration with enough simulations in order to get convergence. This is due to the $B_k$ estimation \eqref{eq:bk} for which we have to take enough simulations to avoid bias propagation.
The result is classical in nested Monte Carlo : not enough convergence in inner iteration can lead to a bias on upper iterations.\\
The convergence result is quite obvious: 
\begin{itemize}
\item
the bias propagates multiplied at each switching dates by a square of the Lipschitz constant but decrease due to the fact that the probability that the branching dates doesn't reach $T$ goes to zero. In fact using  Stirling
formula $ \Gamma(p) \simeq \sqrt{2 \pi (p-1)} \big( \frac{p-1}{e}\big)^{p-1}$ we see that the bias term goes to $0$ exponentially fast meaning that for not too long maturities only small values of $p$ are needed to reach a very good accuracy.
\item The variance term can be bounded by
\begin{align*}
\sum_{i=0}^{p-1} \frac{C_i}{N_i}
\end{align*}
with $C_i$ going to zero very quickly meaning that we should take $N_i$ with decreasing values.\\
Besides if we consider series $\{(N_0^j, .., N_{p-1}^j)\}_{j>0}$ such that the corresponding $(\bar u^p_\emptyset)^j$ goes to the bias term, it is reasonable to take $(N_0^j, .., N_{p-1}^j) = M^j (N_0^0, .., N_{p-1}^0)$ where $M^j$ goes to infinity.
\end{itemize}
\begin{remark}
For $f$ regular  it could be tempting to try to use the  ideas developed in 
\cite{rainforth2017opportunities}.
The approximation in this article for $f$ regular uses the fact the bias goes to $0$ in order to declare that
terms depending on the square of the bias are negligible  compared to terms depending  on the bias at each iteration of the nesting procedure. This is not true in our case. 
\end{remark}
\begin{remark}
The previous result in proposition \ref{theo1} is also valid for more complex SDE as soon as the SDE can be simulated exactly.
\end{remark}
\begin{remark}
We could have a tighter expression for the variance terms by keeping \\
$||g(X_T) -\E\big( g(X_T) \big) ||_2^2$ and $||f(t,X_t,u(t,X_t))-\E\big( f(t,X_t,u(t,X_t))\big)||_2^2$ instead of $||  g(X_T)||_2^2$ and $||f(t,X_t,u(t,X_t))||_2^2$ respectively.
\end{remark}
\begin{remark}
When the coefficients $\mu$ and $\sigma$ are not constant
the methodology is exactly the same except that the SDE has to be
approximated by an Euler Scheme. Two ways to implement it can be used :
\begin{itemize}
\item the first consists in getting the coefficients of the SDE on a fixed grid with a given time step picking the values from the grid. Then the error added due to the discretization is classical \cite{talay1990expansion}.
\item a second numerically more effective consists in  using an Euler scheme between the switching dates essentially meaning that the Euler grid depends on the trajectory. This approach is suggested in \cite{warin2017variations}.
\end{itemize}
\end{remark}
At last we see that the method converges with a speed independent of the dimension of the problem meaning that it is possible to solve non linear PDEs in very high dimension.
\subsection{Numerical results for the first non linear case}
We will first study a first toy case in high dimension then we will move to a realistic test case in finance.
\subsubsection{A first toy example}
In this first case we take 2 maturities $T=1$, $T=2$. 
The SDE coefficients are  $\mu = \frac{\mu_0}{d} \un_d $, $\sigma = \frac{\sigma_0}{\sqrt{d}} \I_d$ with $\mu_0= 0.2$, $\sigma_0= 1$.\\
We take for $x \in \R^d$, $g(x)= \cos(\sum_{i=1}^d x_i)$ and the non linearity:
\begin{align*}
f(t,x, u)= &  \cos(\sum_{i=1}^d x_i)(a+ \frac{\sigma_0^2}{2}) e^{a (T-t)} 
+ \sin(\sum_{i=1}^d x_i) \mu_0 e^{a (T-t)}
- \\ & r \cos(\sum_{i=1}^d x_i)^2 e^{2 a (T-t)}  + r (- e^{a (T-t)} \vee ( u \wedge e^{a (T-t)}))^2
\end{align*}
with $a =0.1$, $r= 0.1$,$d=100$.\\
Equation \eqref{eqPDE} admits the classical solution $u(t,x)= e^{a(T-t)} cos(
\sum_{i=1}^d x_i)$.
The Lipschitz constant associated to $f$ is $ K=2r e^{a T}$ and the solution is Lipschitz in time ($\theta=1$) with a Lipschitz constant $\hat K =a e^{aT}$.\\
Notice that $$\E(g(X_T)^2) \le 1,$$ and $$ \sup_{t \in[0,T]} \E(f(t,X_t,u(t,X_t))^2) \le (a+ \frac{\sigma_0^2}{2} +\mu_0)^2 e^{2a (T-t)}. $$

In tables \ref{tab:coefCase1T0}, \ref{tab:coefCase1T1} we give the different coefficients associated to error expression in proposition \ref{theo1}:
\begin{itemize}
\item Bias $p$ corresponds to term 
\begin{align}
b(p)=\frac{K^{2p} e^{\lambda T} }{\lambda^{p}} T^{2\theta} \hat{K}^2 \frac{T^p}{p\Gamma(p)},
\label{eq:bp}
\end{align}
\item Var $i$ corresponds to the variance term 
\begin{align}
v(i)= K^{2i} \frac{ T^{i} e^{\lambda T}}{\lambda^{i}} \kappa_i.
\label{eq:vi}
\end{align}
\end{itemize}
\begin{table}[h!]
\centering
 \begin{tabular}{|c|c|c|c|}  \hline
 $\lambda $  &0.2 &0.4 &0.8    \\ \hline
 Bias  p=  1 &0.01458 &0.008902 &0.00664    \\ \hline
 Bias  p=  2 &0.00178 &0.0005437 &0.0002028    \\ \hline
 Bias  p=  3 &0.000145 &2.213e-05 &4.128e-06    \\ \hline
 Bias  p=  4 &8.854e-06 &6.759e-07 &6.302e-08    \\ \hline
 Bias  p=  5 &4.326e-07 &1.651e-08 &7.697e-10    \\ \hline
 Var i = 0 &21.54 &14.65 &13.15    \\ \hline
 Var i = 1 &2.929 &1.077 &0.5374    \\ \hline
 Var i = 2 &0.2628 &0.05125 &0.01371    \\ \hline
 Var i = 3 &0.01753 &0.001791 &0.0002515    \\ \hline
 Var i = 4 &0.0009291 &4.93e-05 &3.588e-06    \\ \hline
 \end{tabular}
\caption{\label{tab:coefCase1T0} Coefficient in the error analysis in proposition \ref{theo1} for $T=1$}
\end{table}
 \begin{table}[h!]
\centering
 \begin{tabular}{|c|c|c|c|}  \hline
 $\lambda $  &0.2 &0.4 &0.8    \\ \hline
 Bias  p=  1 &0.2125 &0.1585 &0.1764    \\ \hline
 Bias  p=  2 &0.0634 &0.02364 &0.01316    \\ \hline
 Bias  p=  3 &0.01261 &0.002352 &0.0006542    \\ \hline
 Bias  p=  4 &0.001881 &0.0001754 &2.44e-05    \\ \hline
 Bias  p=  5 &0.0002245 &1.047e-05 &7.28e-07    \\ \hline
 Var i = 0 &59.96 &46.95 &57.2    \\ \hline
 Var i = 1 &18.78 &7.668 &5.005    \\ \hline
 Var i = 2 &3.912 &0.8287 &0.2856    \\ \hline
 Var i = 3 &0.6101 &0.06674 &0.01202    \\ \hline
 Var i = 4 &0.07596 &0.004276 &0.0003996    \\ \hline
 \end{tabular}
\caption{\label{tab:coefCase1T1} Coefficient in the error analysis in proposition \ref{theo1} for $T=2$}
\end{table}

 We check that the bias and the variance  terms  decrease rapidly with $p$ for small maturities.\\
 In tables  \ref{tabNPart1} and \ref{tabNPart2}, for $T=1$ and $T=2$, $\lambda =0.4$, we give the level $p$ and the number of particles to take at each level to reach a given accuracy.

\begin{table}[h!]
\centering
 \begin{tabular}{|c|c|c|}  \hline
 $i$  &0 &1    \\ \hline
 $N_i$  &129684 &5299    \\ \hline
 \end{tabular}
\caption{\label{tabNPart1} Number of particles to take for $ p = 2 ,\lambda=  0.4 ,  T=  1 $, an accuracy  $B_p +   \sum_{i=0}^{p-1} \frac{v(i)}{N_{i}}\prod_{j=0}^{i}(1+ \frac{8}{N_{j-1}})= 1.088E-03$.  }
\end{table}

\begin{table}[h!]
\centering
 \begin{tabular}{|c|c|c|c|}  \hline
 $i$  &0 &1 &2    \\ \hline
 $N_i$  &59885 &9780 &1057       \\ \hline
 \end{tabular}
\caption{\label{tabNPart2} Number of particles to take for $ p = 3 ,\lambda=  0.4 ,  T=  2 $, an accuracy   $B_p + \sum_{i=0}^{p-1} \frac{v(i)}{N_{i}}\prod_{j=0}^{i}(1+ \frac{8}{N_{j-1}})= 4.720E-03$.  }
\end{table}

\begin{figure}[h!]
 \begin{minipage}[b]{0.49\linewidth}
  \centering
 \includegraphics[width=\textwidth]{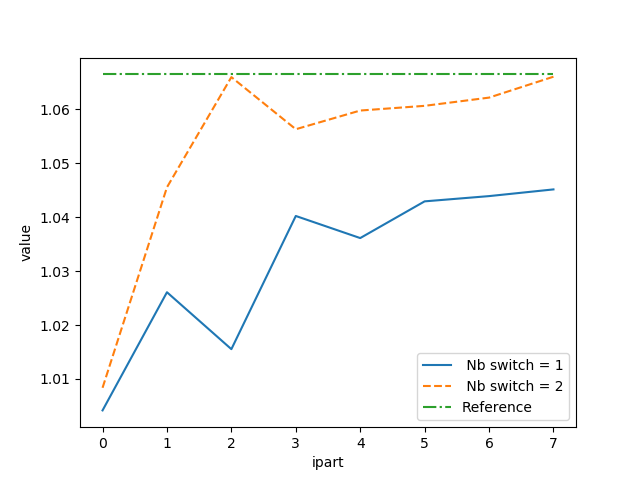}
 \caption*{$\lambda=0.2$.}
 \end{minipage}
\begin{minipage}[b]{0.49\linewidth}
  \centering
 \includegraphics[width=\textwidth]{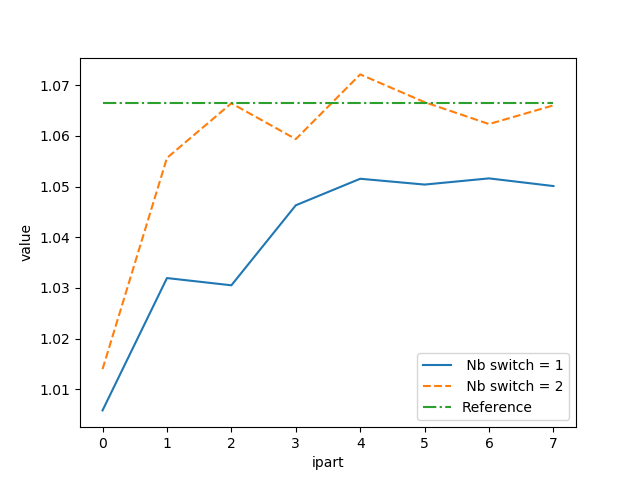}
 \caption*{$\lambda=0.4$.}
 \end{minipage}
 \begin{minipage}[b]{\linewidth}
  \centering
 \includegraphics[width=0.5\textwidth]{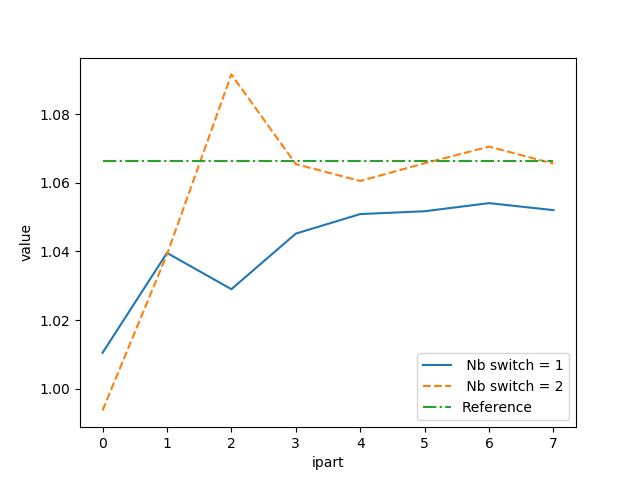}
 \caption*{$\lambda=0.8$.}
 \end{minipage}
\caption{\label{fig:case1T1} Convergence for different number of switches for case 1, $T=1$.}
\end{figure}

On figure \ref{fig:case1T1}, we plot for different values of $\lambda$ the solution  with $T=1$ obtained with  one or two switches with a number of particles $ N_0 = 1000 \times 2^{\mbox{ ipart}}$ and $N_1= 50 \times 2^{\mbox{ ipart}}$. With one switch the solution is clearly biased while the bias is indistinguishable from 0 with $2$ switches whatever the $\lambda$ taken. \\

\begin{figure}[h!]
 \begin{minipage}[b]{0.49\linewidth}
  \centering
 \includegraphics[width=\textwidth]{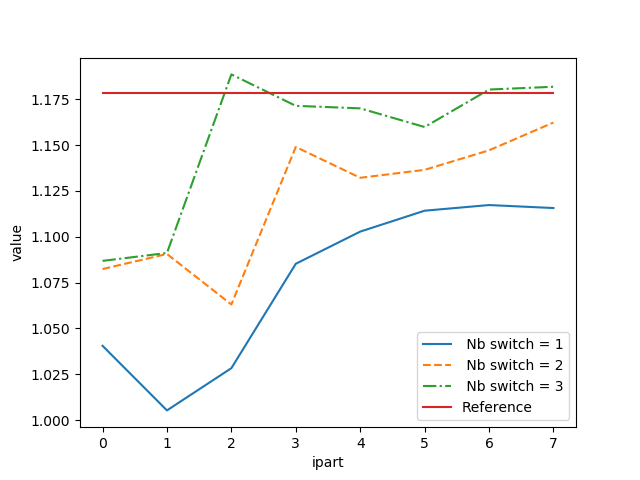}
 \caption*{$\lambda=0.2$.}
 \end{minipage}
\begin{minipage}[b]{0.49\linewidth}
  \centering
 \includegraphics[width=\textwidth]{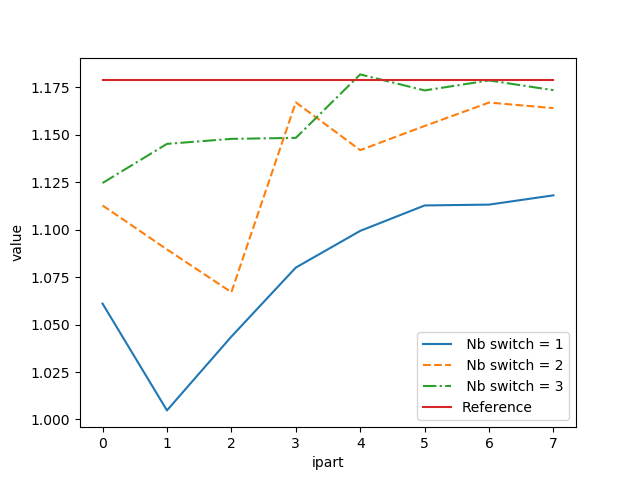}
 \caption*{$\lambda=0.4$.}
 \end{minipage}
 \begin{minipage}[b]{\linewidth}
  \centering
 \includegraphics[width=0.5\textwidth]{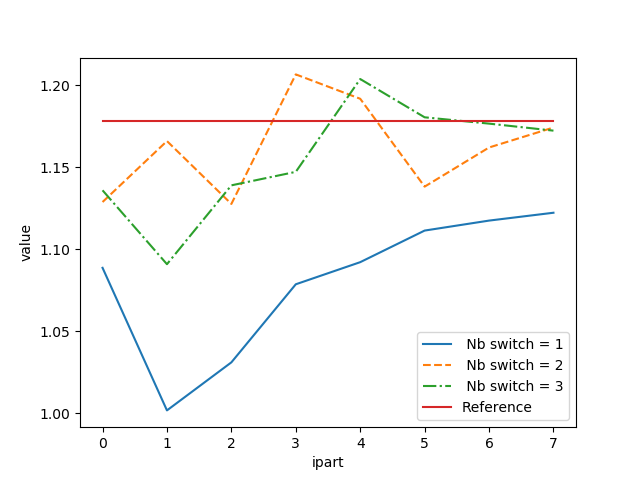}
 \caption*{$\lambda=0.8$.}
 \end{minipage}
\caption{\label{fig:case1T2} Convergence for different number of switches for case 1, $T=2$.}
\end{figure}

On figure \ref{fig:case1T2}, we plot for different values of $\lambda$ the solution with $T=2$ obtained with  one, two or three switches with a number of particles $ N_0 = 1100 \times 2^{\mbox{ ipart}}$, $N_1= 110 \times 2^{\mbox{ ipart}}$, $N_2=25 \times 2^{\mbox{ ipart}} $. For all $\lambda$, we have to take three switches to have a good precision. The best solution seem to be reached with $\lambda=0.2$.
\newpage
\subsubsection{A second test case}
\begin{figure}[h!]
\begin{minipage}[b]{0.49\linewidth}
  \centering
 \includegraphics[width=\textwidth]{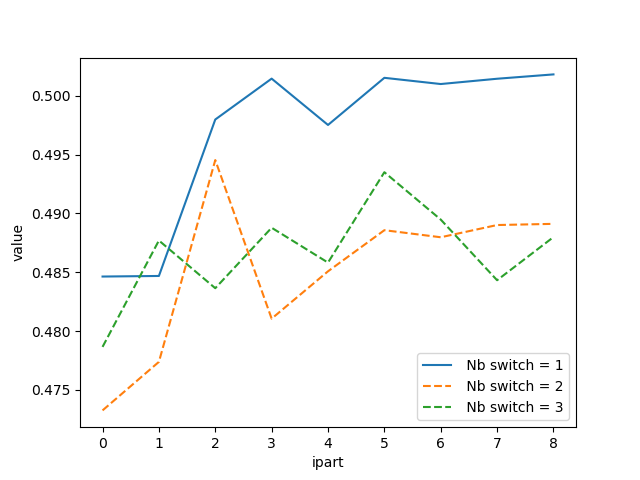}
 \caption*{$\lambda=0.1$.}
 \end{minipage}
\begin{minipage}[b]{0.49\linewidth}
  \centering
 \includegraphics[width=\textwidth]{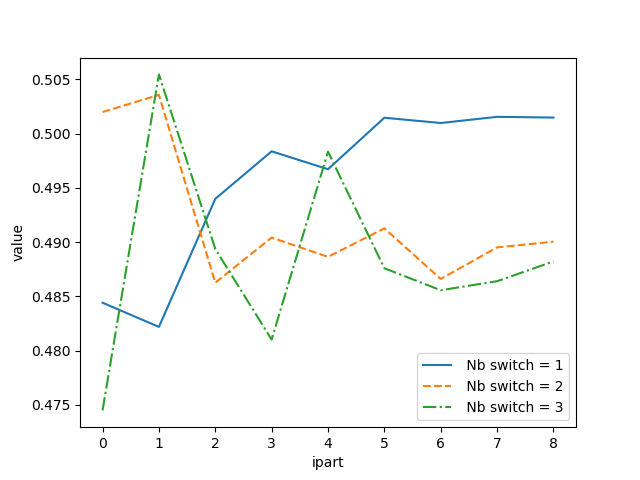}
 \caption*{$\lambda=0.2$.}
 \end{minipage}
\caption{\label{fig:caseLabor} Convergence of the scheme on the CVA case.}
\end{figure}

We use the test case presented in \cite{henry2017deep}.
This is a test case in low dimension but the author gives some numerical bounds on the solution so that we can compare our methodology to some deep learning solution.\\
The author considers the PDE obtained from a CVA valuation problem.
 \begin{flalign}
 \Lc u(t,x) := \mu Du(t,x) +  \frac{1}{2} \sigma \sigma^{\top} \!:\! D^2 u(t,x),
 \label{eq:labor}
 \end{flalign}
 taking $\mu= -\frac{\sigma_0^2}{2} \un_d$, $\sigma = \sigma_0 \I_d$, and a non-linearity
 \begin{align*}
 f(t,x,u)=  \beta (u^+ -u),
 \end{align*}
 with $\beta=0.03, \sigma_0 =0.2$.
 The initial value for the SDE is  $X_0=  \un_d$ .
 The final function $g(x)=\sum_{i=1}^d (1-2 1_{e^{x_i} >1})$  and $T=1.$
 Some bounds on the solution in dimension till 6 are given. For $d=6$ a lower bound calculated is $48.80$, an upper bound is $48.83$ whereas with $\beta=0$ the solution is $47.73$.\\
 On this simple problem we don't try to optimize  the number of particles taken at each level nor the $\lambda$ taken.
 On figure \ref{fig:caseLabor}, we plot for $\lambda=0.1$ and $\lambda=0.2$, the solution  obtained with  one, two  or three switches with a number of particles $ N_0 = 36000 \times 2^{\mbox{ ipart}}$, $N_1= 140 \times 2^{\mbox{ ipart}}$, $N_2= 2^{\mbox{ ipart}}$.
With  $\lambda =0.1$ with 3 switches and $ipart=8$ we get $0.4880$ while with $\lambda=0.2$ we get $0.4882$ so that both values are in the very tight bounds proposed in \cite{henry2017deep}.

\subsubsection{A third test case}
\begin{figure}[h!]
  \centering
 \includegraphics[width= 0.5 \textwidth]{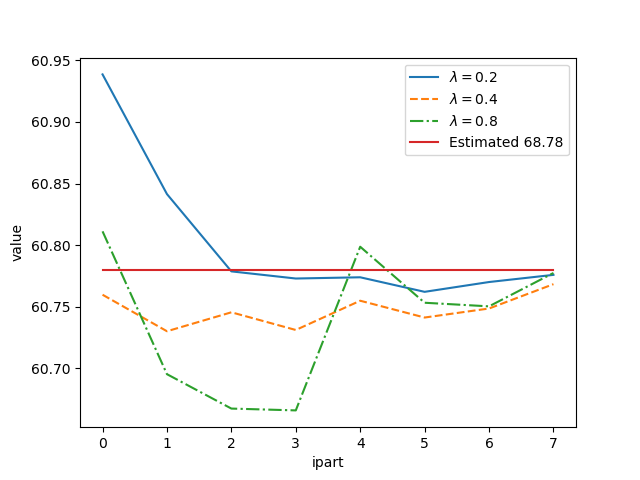}
\caption{\label{fig:case2Nlin} Convergence of the scheme on the Black Scholes case without default  taking $f=0$ in the algorithm with a number of particle $N = 100000 \times 2^{\mbox{ipart}}$.}
\end{figure}

In this part we take a test case coming from \cite{han2017overcoming} modeling the valuation of an European claim in dimension 100 using a  Black-Scholes dynamic of the assets supposing the existence of  a default risk.
The default is modeled  by the first jump time of a Poisson process with intensity Q.
When a default occurs, the claim's holder receive only a fraction $\delta \in [0,1]$ of the current value.
We want to valuate  the claim conditionally that the default hasn't occurred yet.
The dynamic of an asset $S_t$ with trend $\mu_0$ and volatility $\sigma_0$ following the BS model  satisfies
\begin{align}
\label{eq:BS}
S_t = S_0 e^{(\mu_0 - \frac{\sigma_0^2}{2}) t + \sigma_0 W_t},
\end{align}
such that taking $X_t = \log(S_t)$, the $X_t$ dynamic follows
\begin{flalign}
 d X_t =  (\mu_0- \frac{\sigma_0^2}{2}) dt+ \sigma_0 dW_t.
 \end{flalign}
 Supposing that all the assets are independent and follow the same equation \eqref{eq:BS},
the value of the claim given by \cite{han2017overcoming} can be equivalently given as the solution at date $0$ and point $x = \log(100)\un_d$ of \eqref{eqPDE} where
 \begin{flalign}
 \Lc u(t,x) := \mu Du(t,x) +  \frac{1}{2} \sigma \sigma^{\top} \!:\! D^2 u(t,x),
 \label{eq:genBS}
 \end{flalign}
 with $\mu = (\mu_0-\frac{\sigma_0^2}{2}) \un_d $, $\sigma = \sigma_0 \I_d$.\\
 Following \cite{han2017overcoming}, the final function $g$ satisfies for $x \in \R^d$, $$g(x)=\min_{i=1}^{100}(e^{x_i}),$$
 while the non-linearity is given by
 \begin{align*}
 f(t,x,u)= - \left( (1-\delta ) \min\{ \gamma^h, \max\{ \gamma^l,\frac{\gamma^h-\gamma^l}{v^h-v^l}(u-v^h) + \gamma^h\}\} +R \right) u
 \end{align*}
 where $R$ is the interest rate of the riskless asset (see \cite{bergman1995option}).
 We take the same parameters as in \cite{han2017overcoming} so $T=1$, $\delta =\frac{2}{3}$, $\mu_0=0.02$, $\sigma_0=0.2$, $v^h=50$, $v^l=70$, $\gamma^h=0.2$, $\gamma^l=0.02$.
 As pointed out by \cite{han2017overcoming}, the solution obtained by Monte Carlo ignoring the default risk is approximately $60.78$. This reference can also be obtained by the algorithm taking $f =0$ as shown on figure
 \ref{fig:case2Nlin}.

Let us study the value of the different terms :
\begin{itemize}
\item $ K = (1-\delta) \gamma^h +R =0.086$,
\item $\E(g(X_T)^2) \le  \E((X^1_T)^2) = S_0^2 e^{(2\mu+\sigma_0^2) T} \simeq  10000$  where $X^1$ stands for the first asset log value,
\item Noting $\hat u$ the solution without default
\begin{align*}
\sup_{t \in [0,1]} \E(f(t,X_t, u(t,X_t))^2) \le &  K^2 \sup_{t \in [0,1]} \E(  u(t,X_t)^2)  \\
\le & K^2 \sup_{t \in [0,1]} \E( \hat u(t,X_t)^2) \le K^2 \sup_{t \in [0,1]} \E((X^1_t)^2) \\
= &K^2 S_0^2 \sup_{t \in [0,1]} e^{(2\mu+\sigma_0^2) t}\simeq 50.
\end{align*}
\item At last the solution may not be uniformly Lipschitz in time but remember that
$ \hat K^2 T^{2 \theta} $ is a bound from an expression  $\psi = \sup_{t \in [0,T]} \E( (\hat u(t,X_t) -g(X_t))^2)$ so that using the previous estimations 
\begin{align*}
\psi \le 2 \sup_{t \in [0,T]} \E( u(t,X_t)^2  + g(X_t)^2) \\
\le 2 S_0^2 \sup_{t \in [0,1]} (e^{(2\mu+\sigma_0^2) t} +e^{(2\mu+\sigma_0^2) t} ) \simeq  40000
\end{align*}
so that $ \hat K^2 T^{2 \theta} $ can be replace $40000.$
\end{itemize}
In table \ref{tab:coefCase2}, we give the coefficients of equations \eqref{eq:bp} and \eqref{eq:vi} involved in proposition \ref{theo1}.
\begin{table}[h!]
\centering
\begin{tabular}{|c|c|c|c|}  \hline
 $\lambda $  &0.2 &0.4 &0.8    \\ \hline
 Bias  p=  1 &1988.0 &1214.0 &905.4    \\ \hline
 Bias  p=  2 &37.32 &11.4 &4.25    \\ \hline
 Bias  p=  3 &0.4672 &0.07134 &0.0133    \\ \hline
 Bias  p=  4 &0.004387 &0.0003349 &3.122e-05    \\ \hline
 Bias  p=  5 &3.295e-05 &1.258e-06 &5.863e-08    \\ \hline
 Var i = 0 &28450.0 &33540.0 &49120.0    \\ \hline
 Var i = 1 &1031.0 &618.3 &457.0    \\ \hline
 Var i = 2 &19.13 &5.77 &2.139    \\ \hline
 Var i = 3 &0.238 &0.036 &0.006682    \\ \hline
 Var i = 4 &0.002226 &0.0001687 &1.567e-05    \\ \hline
 \end{tabular}
 \caption{ \label{tab:coefCase2}
 Coefficient in the error analysis in proposition \ref{theo1} for the Black Scholes case with default risk.}
\end{table}
By taking $\lambda =0.8$ , $p=3$, on table \ref{tabNPartCase2} we give the number of particles to take  to have an accuracy of $0.01$.
\begin{table}[h!]
\centering
 \begin{tabular}{|c|c|c|c|}  \hline
 $i$  &0 &1 &2    \\ \hline
 $N_i$  &  110796 &1030 &4  \\ \hline
 \end{tabular}
\caption{\label{tabNPartCase2} Number of particles to take for $ p = 3 ,\lambda=  0.8 $, an accuracy  $B_p + \sum_{i=0}^{p-1} \frac{v(i)}{N_{i}}\prod_{j=0}^{i}(1+ \frac{8}{N_{j-1}})= 4.1E-2. $  }
\end{table}

\begin{figure}[h!]
\begin{minipage}[b]{0.49\linewidth}
  \centering
 \includegraphics[width=\textwidth]{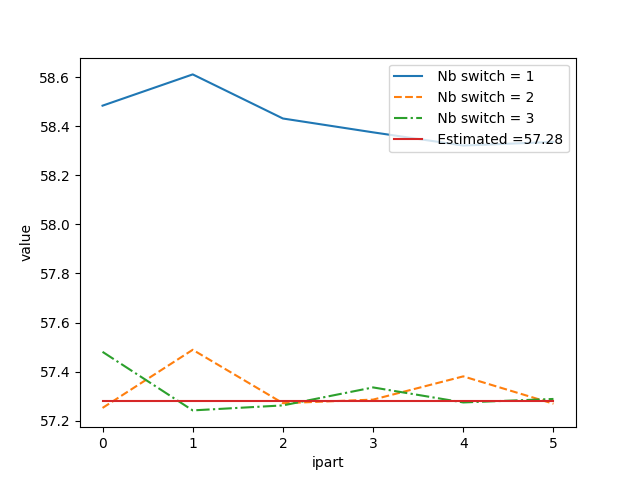}
 \caption*{$\lambda=0.2$.}
 \end{minipage}
\begin{minipage}[b]{0.5\linewidth}
  \centering
 \includegraphics[width=\textwidth]{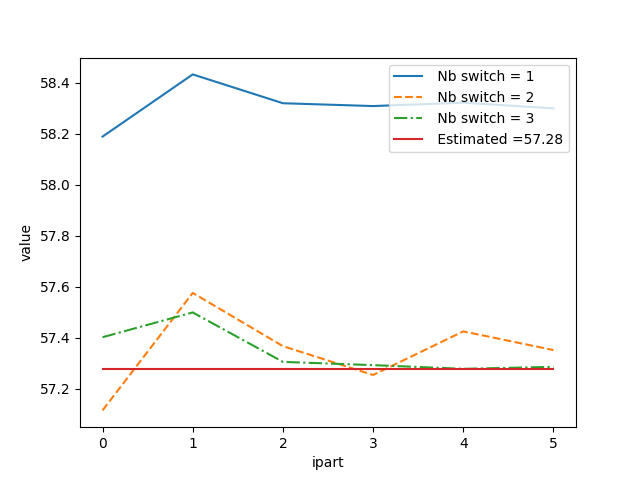}
 \caption*{$\lambda=0.4$.}
 \end{minipage}
 \begin{minipage}[b]{\linewidth}
  \centering
 \includegraphics[width=0.5\textwidth]{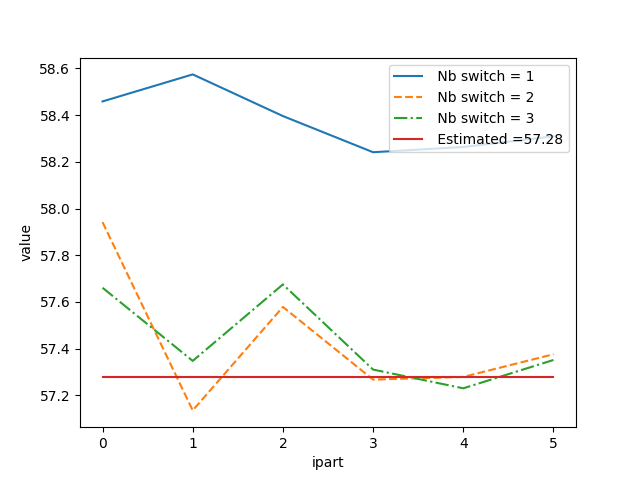}
 \caption*{$\lambda=0.8$.}
 \end{minipage}
\caption{\label{fig:case2} Convergence of the scheme on the Black Scholes case with default}
\end{figure}
On figure \ref{fig:case2}, we plot for different values of $\lambda$ the solution  obtained with  one, two  or three switches with a number of particles $ N_0 = 36000 \times 2^{\mbox{ ipart}}$, $N_1= 40 \times 2^{\mbox{ ipart}}$, $N_2= 2^{\mbox{ ipart}}$.
The solution seems to be $57.28$  (value obtained for $3$ switches with both $\lambda=0.1$ and $\lambda=0.2$  and also obtained with deep learning techniques \cite{Quentin}) and close to the value obtained in \cite{han2017overcoming} who game $57.30$.\\
Two switches are enough to get a good accuracy. For example 
$57.27$ is reached  with two switches taking $N_0=1152000$, $N_1= 4480$ in $90$ seconds with $\lambda=0.2$.
\newpage
\section{The semi linear case}
In this section we extend the previous scheme obtained to the semi-linear case.
To simplify the setting, without restriction, we just take a function depending on $Du$:
 \begin{flalign}
 \label{eqPDESemi}
  (-\partial_tu-\Lc u)(t,x)  & = f(t,x,Du(t,x)), \nonumber \\
  u_T&=g, \quad  t<T,~x\in\R^d,
 \end{flalign}
 where $\Lc$ is always given by equation \eqref{eq:gen} and 
 the dynamic of the underlying  SDE is still given by \eqref{eq:sde} 
 with $\mu \in \R^d$, and $\sigma \in \M^d$ is here some constant non-degenerated matrix.\\
We will take the same kind of assumption as in the previous section:
 \begin{assumption}
 \label{ass::lipFSemi}
 $f$ is  uniformly Lipschitz in $Du$ with constant $K$ :
 \begin{flalign}
 | f(t, x, y) - f(t, x, w) | \le K ||y-w||_2 \quad \quad \forall t \in [0,T], x \in \R^d, (w,y) \in \R^d \times \R^d.
 \end{flalign}
 \end{assumption}

 \begin{assumption}
 \label{ass::uHolderSemi}
 Equation  \eqref{eqPDESemi} has a solution $u \in C^{1,2}([0,T] \x \R^d)$, such that
 \begin{itemize}
 \item $Du$ is $\theta$-H\"older with $\theta \in (0,1]$ in time with constant $\hat K$ :
 \begin{align*}
 || Du(t,x) - Du(\tilde t, x) || \le \hat K | t- \tilde{t}|^\theta \quad \quad \quad \forall (t,\tilde{t},x) \in [0,T] \times [0,T] \times \R^d,
 \end{align*}
 \item $u(t,x)$ and $Du(t,x)$ have a quadratic growth in $x$ uniformly in $t$.
 %\item $\E\big( \int_0^T |f(t,X_t, Du(t,X_t))| dt\big)< \infty$.
 \end{itemize}
 \end{assumption}
 \begin{assumption}
 \label{ass::gLipSemi}
 $g$ is uniformly Lipschitz  such that for $\tilde K >0$:
 \begin{align*}
 |g(x) -g(y)| \le \tilde K ||x-y||_2 \quad \forall (x,y) \in \R^d \times \R^d.
 \end{align*}
 \end{assumption}
 \subsection{General idea of the algorithm}
 We will propose two algorithms that are some extensions of the algorithm previously given.\\
 As in the previous section, the sequence $(T_i)_{i \ge 0}$ is defined by equation \eqref{eq:T} but the $(\tau_{m})_{m \ge 1}$ are i.i.d. random variables of density $\rho$ which follow  a general  gamma distribution so that
 \begin{align}
 \rho(x)= \lambda^u x^{u-1} \frac{e^{-\lambda x}}{\Gamma(u)}, u >0
 \label{rho}
 \end{align}
 and the associated  cumulated distribution function is   $$F(x) =\frac{\gamma(u,\lambda x)}{\Gamma(u)},$$
 where $\gamma(s, x) =  \int_0^x t^{s-1} e^{-t} dt$ is the incomplete gamma function.\\
  In order to have a converging method we will see that we will have to take $u<1$ in $\rho$ expression \eqref{rho} excluding the exponential distribution. This a weaker constraint than in \cite{henry2016branching} where, using branching for some  polynomial non-linearities, converging results were only obtained for $u<0.5$.
 \\
 
 Under the regularity assumption on $u$, from the Feynman-Kac formula,  the representation of the solution $u$ is
\begin{flalign}
 u(0,x) 
 = & 
 \E_{0,x} \Big[ \Fb(T)\frac{g(W_T)}{ \Fb(T)}+\int_0^T \frac{f(t,X_t,Du(t,X_t))}{\rho(t)}\rho(t)dt\Big] \nonumber 
 \\ =&
 \E_{0,x} \big[ \hat \phi\big(0, T_{1},X_{T_{1}}, Du(T_{1},X_{T_{1}  })\big)\big],
 \label{eq:recSemi}
 \end{flalign}
 with
  \begin{flalign}
  \hat \phi(s, t,x,z) &:= \frac{\1_{\{t\ge T\}}}{ \Fb(T-s)} g(x)\!+\! \frac{\1_{\{t<T\}}}{\rho(t -s)} f(t,x,z).
 \label{eq:phihat}
 \end{flalign}
 then we define  $Du(T_{1},X_{T_{1}})$ using the automatic differentiation rule :
 \begin{align}
 Du(T_{1},X_{T_{1}}) = \E_{T_{1},X_{T_{1}}} \Big[ \sigma^{-\top} \frac{W_{T_{2}}-W_{T_{1}}}{T_{2}-T_{1}} \phi(T_{1}, T_{2},X_{T_{1}},X_{T_{2}},Du(T_{2},X_{T_{2}}))\big],
 \label{eq:ruelDu}
 \end{align}
 with
  \begin{flalign}
  \phi(s, t,x,y,z) &:= \frac{\1_{\{t\ge T\}}}{ \Fb(T-s)}(g(y)-g(x))\!+\! \frac{\1_{\{t<T\}}}{\rho(t -s)} f(t,y,z).
 \label{phi}
 \end{flalign}
 where  the $g(x)$ acts as a control variate term.\\
 The automatic differentiation used here is 
based on the Malliavin integration by parts formula  (see \cite{fournie1999applications} for its use
in the context of Monte Carlo approximation and the extension to other sensitivities) and has been used in a similar context as the one presented here  in \cite{doumbia2017unbiased}, \cite{henry2016branching}.\\
  Recursively we define for $n < N_T$:
 \begin{flalign}
 Du_{n} = & \E_{ T_{n}, X_{T_n}}\big[ \sigma^{-\top} \frac{W_{T_{n+1}}-W_{T_{n}}}{T_{n+1}-T_{n}} \phi\big(T_n , T_{n+1}, X_{T_{n}}, X_{T_{n+1}}, Du_{n+1}\big)\big],
 \end{flalign}
 As in the previous section we besides consider the truncated operator after $p$ switches:
 \begin{flalign*}
 u_0^p =& \E \big( \hat \phi\big(0, T_{1}, X_{T_{1}}, Du_{1}^p\big) \big)\\
  Du_{n}^p = & \E_{ T_{n}, X^n_{T_n}}\big[ \sigma^{-\top}  \frac{W_{T_{n+1}}-W_{T_{n}}}{T_{n+1}-T_{n}}\phi\big(T_n, T_{n+1}, X_{T_{n}}, X_{T_{n+1}}, Du_{n+1}^p \big)\big], \quad 1 \le n < p \\
  Du_p^p = & Dg(X_{T_p}).\\
 \end{flalign*}
 The goal of the following section is to present two algorithms based on the previously defined recursion and to show their convergence.
 
 \subsection{A first estimator}
 \label{sec:fistEstim}
 We take the same notations as in the  section \ref{sec::estim} for the set $Q_i$, $i<p$, the set $\tilde Q(k)$  for $k \in Q_i$.
 The $\tau_{k}$ are as before some switching increments. They are  always i.i.d. random variables with density $\rho$ and for $k\in Q^p$  the $\bar W^{k}$ are some independent $d$-dimensional Brownian motions,  independent of the $(\tau_{k})_{k \in Q^p}$ too.
 The switching dates are defined  by equation \eqref{eq:T} and $(X^{\tilde k}_t)_{t\geq 0}$ defined by
 \eqref{eq:diff}. \\
 We propose the following estimator
 defined  by:
\begin{align}
\left \{
\begin{array}{ll}
\bar u_\emptyset^p = & \frac{1}{N_0} \sum_{j=1}^{N_0} \hat \phi\big( 0,T_{(j)}, X^{(j)}_{T_{(j)}}, D \bar  u_{(j)}^p\big), \\
 D \bar u_{k}^p = &  \frac{1}{N_i} \sum_{\tilde k \in \tilde Q(k)} \phi\big(T_k, T_{\tilde k}, X^{k}_{T_{k}}, X^{\tilde k}_{T_{\tilde k}},  D\bar u_{\tilde k}^p\big)
\sigma^{-\top} \frac{\bar W_{T_{\tilde{ k}}- T_k}^{\tilde k}}{T_{\tilde{ k}}- T_k} \\
& \quad \mbox{ for } k =( k_1, k_2,..k_i) \in Q_i, i <p, \\
D \bar u_{\tilde k}^p = &   Dg(X^{\tilde k}_{T_{\tilde k}}) \quad \mbox{ for } 
\tilde k \in Q_p
\end{array}
\right.
\label{eq:estimSemi1}
\end{align}
\begin{remark}
An estimator of the gradient at the initial date is off course  available too as
\begin{flalign*}
D \bar u_{\emptyset}^p = &  \frac{1}{N_0} \sum_{j=1}^{N_0}  \phi\big(0,T_{(j)}, x, X^{(j)}_{T_{(j)}}, D \bar u_{(j)}^p\big)
\sigma^{-\top} \frac{\bar W_{T_{(j)}}}{T_{(j)}}
\end{flalign*}
\end{remark}

\begin{proposition}
\label{theo2}
Suppose that $\rho$ is the density of a gamma law so that $ \rho(x)=\lambda^u x^{u-1} \frac{e^{-\lambda x}}{\Gamma(u)}$,  suppose that  $u<1$, and suppose that assumptions \ref{ass::lipFSemi},  \ref{ass::uHolderSemi},  \ref{ass::gLipSemi}   are satisfied then  there are two  functions $C$ and $\tilde C$ depending on $\sigma$ and one $\bar C$ depending on $\mu$, $\sigma$ and $T$ such that using estimator \eqref{eq:estimSemi1} to solve equation  \eqref{eqPDESemi}, we  have the error estimate
\begin{align}
\label{est::finalSemi1}
 \E\big((\bar u_{\emptyset}^p - u(0,x))^2 \big) \le
 \prod_{i=1}^p (1+\frac{8}{N_{i-1}}) \frac{\Gamma(u)^pe^{\lambda T}}{\lambda^{p}} \frac{T^{(1-u)p+1+2\theta}}{(1-u)^{p-1}(2-u)}C(\sigma)^{p-1} \hat K^2  K^{2p} + \nonumber \\
  4 \sum_{i=0}^{p-1} \frac{K^{2i}}{N_{i}} \prod_{j=1}^i (1+\frac{8}{N_{j-1}})
 \frac{\Gamma(u)^{i+1} e^{\lambda T}}{\lambda^{i+1}} \frac{T^{(1-u)(i+1)+1}}{(1-u)^{i}(2-u)}  \tilde C(\sigma)^{i} \hat F +  \nonumber  \\
 2 \sum_{i=1}^{p-1} \frac{K^{2i}}{N_{i}} \prod_{j=1}^i (1+\frac{8}{N_{j-1}})
 \frac{\Gamma(u)^{i+1} e^{\lambda T}}{\lambda^{i}} \frac{T^{(1-u)i+1}}{(1-u)^{i-1}(2-u)} \frac{ \bar C(\mu, \sigma,T)  C(\sigma)^{i-1} \tilde K^2}{ \Gamma(u)-\gamma(u,\lambda T )}  + \nonumber \\
  \frac{2}{ N_0} \frac{\Gamma(u)}{\Gamma(u)-\gamma(u,\lambda T)}   \E(g(X_T)^2)
\end{align}
where
\begin{align}
  \label{eq:FHat}
  \hat F = \sup_{t \in [0,T]} \E[ f(t,X_t,Du(t,X_t))^4]^{\frac{1}{2}},
\end{align}
\end{proposition}
\begin{proof}
Under assumption \ref{ass::uHolderSemi}, the solution $u$ of \eqref{eqPDESemi} satisfies  a Feynman-Kac relation   so that for all $\tilde{k} \in \tilde{Q}(\emptyset)$ 
\begin{align}
u(0,x) =& \E_{ 0, x} \big[  \hat \phi\big(0, T_{\tilde k}, X_{T_{\tilde k}}^{\tilde{k}}, Du(T_{\tilde k},X_{T_{\tilde k}}^{\tilde{k}} ))\big],
\label{eq:realU0}
\end{align}
and that for all $k \in Q_i$ with $ i< p$,  and $\forall \tilde k \in \tilde{Q}(k)$, the gradient is given by:
\begin{align}
Du(T_k, X^k_{T_k}) =& \E_{ T_{k}, X^k_{T_k}} \big[ \sigma^{-\top} \frac{\bar W^{\tilde{k}}_{T_{\tilde{ k}}- T_k}}{T_{\tilde{ k}}- T_k} \phi\big(T_k , T_{\tilde k},X_{T_{k}}^{k}, X_{T_{\tilde k}}^{\tilde{k}},  Du(T_{\tilde k},X_{T_{\tilde k}}^{\tilde{k}} ))\big],
\label{eq:realU}
\end{align}

We then introduce for $ k \in Q_i$, $1 \le i < p$:
\begin{flalign*}
E_k := & \E_{T_{k},X^{k}_{T_{k}} }\big( || D \bar u_{k}^p - Du(T_k,X^k_{T_k})||^2_2  1_{T_k <T} \big).
\end{flalign*}
In exactly the same way as in the demonstration of proposition \ref{theo1}, we have the following result similar to the one given by equation \eqref{eq:ErrLess}:
\begin{align}
E_k \le& \frac{1}{N_i} (1+\frac{8}{N_i}) \sum_{\tilde k \in \tilde Q(k)} E_{T_k,X^{k}_{T_{k}}} \big( \frac{(\bar W^{\tilde{k}}_{\tau_{\tilde{k}}})^\top \sigma^{-1} \sigma^{-\top} \bar W^{\tilde{k}}_{\tau_{\tilde{k}}}}{\tau_{\tilde{k}}^2} \frac{K^2}{\rho(\tau_{\tilde k})^2} E_{\tilde k} \big) + \nonumber
\\
& 4 \frac{1}{N_i^2} \sum_{\tilde k \in \tilde Q(k)} \E_{T_k,X^k_{T_k} }\big(  1_{T_{\tilde k < T}} \frac{(\bar W^{\tilde{k}}_{\tau_{\tilde{k}}})^\top \sigma^{-1} \sigma^{-\top} \bar W^{\tilde{k}}_{\tau_{\tilde{k}}}}{\tau_{\tilde{k}}^2} (\frac{f(T_{\tilde k}, X^{\tilde k}_{T_{\tilde k}}, Du(T_{\tilde k},X^{\tilde k}_{T_{\tilde k}}))}{\rho(\tau_{\tilde k})} )^2\big) +  \nonumber \\
&2 \frac{1}{N_i^2} \sum_{\tilde k \in \tilde Q(k)} \E_{T_k,X^k_{T_k} }\big( 1_{T_{\tilde k\ge T}} \frac{(\bar W^{\tilde{k}}_{T-T_k})^\top \sigma^{-1} \sigma^{-\top} \bar W^{\tilde{k}}_{T-T_k}}{(T-T_k)^2} \frac{(g(X^{\tilde k}_T) -g(X^k_{T_k}))^2}{ \Fb(T-T_k)^2} \big)
\label{eq:ErrLessD}
\end{align}
Using equation \eqref{eq:ErrLess} with $k= \emptyset$ obtained in demonstration of proposition \ref{theo1}:
\begin{align}
 \E\big((\bar u_{\emptyset}^p - u(0,x))^2 \big) \le& \frac{1}{N_0} (1+\frac{8}{N_0}) \sum_{\tilde k \in \tilde Q(\emptyset)} \E \big(  \frac{K^2}{\rho(\tau_{\tilde k})^2} E_{\tilde k} \big) + \nonumber
\\
& 4 \frac{1}{N_0^2} \sum_{\tilde k \in \tilde Q(\emptyset)} \E \big(  1_{T_{\tilde k < T}} (\frac{f(T_{\tilde k}, X^{\tilde k}_{T_{\tilde k}}, Du(T_{\tilde k},X^{\tilde k}_{T_{\tilde k}}))}{\rho(\tau_{\tilde k})} )^2\big) +  \nonumber \\
&2 \frac{1}{N_0^2} \sum_{\tilde k \in \tilde Q(\emptyset)} \E_{T_k,X^k_{T_k} }\big( 1_{T_{\tilde k\ge T}} \frac{g(X^{\tilde k}_T)^2}{ \Fb(T-T_k)^2} \big).
\label{eq:E0}
\end{align}
Iterating  from 1 to $p-1$ we get:
\begin{align}
\label{est::finalErroSemi}
 \E\big((\bar u_{\emptyset}^p - u(0,x))^2 \big) \le& A_1 + 4 A_2 + 2 A_3,
\end{align}
where 
\begin{align}
A_1 = \prod_{i=1}^p \frac{1}{N_{i-1}} (1+\frac{8}{N_{i-1}}) \sum_{\tilde k^1 \in \tilde Q(\emptyset)} ... \sum_{\tilde k^p \in \tilde Q(\tilde{k}^{p-1})}  B_1(\tilde{k}^1,..,\tilde{k}^P)  
\label{eq:A1}
\end{align}
\begin{align}
A_2 = \sum_{i=0}^{p-1} \frac{K^{2i}}{N_{i}^2} \prod_{j=1}^i \frac{1}{N_{j-1}} (1+\frac{8}{N_{j-1}}),
\sum_{\tilde k^1 \in \tilde Q(\emptyset)} ... \sum_{\tilde k^{i+1} \in \tilde Q(\tilde{k}^{i})}  B_2(\tilde{k}^1,..,\tilde{k}^{i+1}),
\label{eq:A2}
\end{align}
and 
\begin{align}
A_3 =\sum_{i=0}^{p-1} \frac{K^{2i}}{N_{i}^2} \prod_{j=1}^i \frac{1}{N_{j-1}} (1+\frac{8}{N_{j-1}})
\sum_{\tilde k^1 \in \tilde Q(\emptyset)} ... \sum_{\tilde k^{i+1} \in \tilde Q(\tilde{k}^{i})}  B_3(\tilde{k}^1,..,\tilde{k}^{i+1}),
\label{eq:A3}
\end{align}
where noting 
\begin{align}
\psi_{k} = \frac{(\bar W^{k}_{\tau_{k}})^\top \sigma^{-1} \sigma^{-\top} \bar W^{k }_{\tau_{k}}}{\tau_{k}^2} , \nonumber \\ 
\hat f_k = f(T_{k}, X^{k}_{T_{k}}, Du(T_{k},X^{k}_{T_{k}})),
\label{eq:not1}
\end{align}
we have:
\begin{flalign}
B_1(\tilde{k}^1,..,\tilde{k}^P) = & \E[   \prod_{j=1}^p \frac{K^2}{\rho(\tau_{\tilde{k}^j})^2} \prod_{j=2}^p \psi_{\tilde{k}^j} E_{\tilde{k}^p}], \nonumber \\
B_2(\tilde{k}^1,..,\tilde{k}^{i+1}) = & \E \big [ 1_{T_{\tilde{k}^{i+1} < T}} (\hat f_{\tilde k^{i+1}})^2
 \prod_{j=1}^{i+1} \frac{1}{\rho(\tau_{\tilde{k}^j})^2} \prod_{j=2}^{i+1} \psi_{\tilde{k}^j} \Big], \nonumber\\
 B_3(\tilde{k}^1,..,\tilde{k}^{i+1}) = &\E \big [1_{T_{\tilde{k}^{i+1} \ge T}} 1_{T_{\tilde{k}^{i} < T}} \frac{(g(X^{\tilde k ^{i+1}}_{T})- 1_{i\ge1} g(X^{\tilde k ^{i}}_{T_{\tilde{k}^i}}))^2}{\bar F(T-T_{\tilde{k}^{i}})^2} \prod_{j=1}^{i} \frac{1}{\rho(\tau_{\tilde{k}^j})^2} \prod_{j=2}^{i} \psi_{\tilde{k}^j} \nonumber \\
 & \frac{(\bar W^{\tilde k^{i+1}}_{T-T_{\tilde k^i}})^\top \sigma^{-1} \sigma^{-\top} \bar W^{\tilde k^{i+1} }_{T-T_{\tilde k^i}}}
 {(T-T_{\tilde k^i})^2}  \big ],
 \label{eq:alllB}
 \end{flalign}
 where $ E_{\tilde{k}^p} = 1_{T_{\tilde k^p < T}}  (Dg(X_{T_{\tilde{k}^p}}^{\tilde{k}^p})- Du(T_{\tilde{k}^p},X_{T_{\tilde{k}^p}}^{\tilde{k}^p}))^2$.\\
We first bound $B_1$.\\
We introduce
$\psi_k = \frac{\phi_k}{\tau_{k}}$,
where 
\begin{align}
\phi_k =G_{k}^\top \sigma^{-1} \sigma^{-\top} G_{k},
\label{eq::phi}
\end{align}
and 
 $G_k \in \R^d$ is composed of centered unitary independent Gaussian random variables.\\
 We introduce $C(\sigma)= \E(\phi_{k})$ which is independent of $k$.
Using the H\"{o}lder property of $Du$, the tower property, the independence of the $\tau_k$ an the $\phi_k$, we get:
\begin{align*}
B_1(\tilde{k}^1,..,\tilde{k}^P)  = & \E[   \prod_{j=1}^p \frac{K^2}{\rho(\tau_{\tilde{k}^j})^2} (\prod_{j=2}^p \frac{1}{\tau_{\tilde{k}^j}}) \E( E_{\tilde{k}^p} \prod_{j=2}^p \phi_{\tilde{k}^j}| \tau_{\tilde{k}^1},.., \tau_{\tilde{k}^p})] \\
\le  &\hat K^2 T^{2\theta} K^{2p} C(\sigma)^{p-1} E[ 1_{T_{\tilde{k}^p < T}}   \prod_{j=1}^p \frac{1}{\rho(\tau_{\tilde{k}^j})^2} \prod_{j=2}^p \frac{1}{\tau_{\tilde{k}^j}} ].
\end{align*}
Using the expression for the density of $\rho$ given by \eqref{rho}, we have the bound:
\begin{align}
H:= &E[ 1_{T_{\tilde{k}^p < T}}   \prod_{j=1}^p \frac{1}{\rho(\tau_{\tilde{k}^j})^2} \prod_{j=2}^p \frac{1}{\tau_{\tilde{k}^j}} ] \nonumber \\
\le  & \frac{\Gamma(u)^p e^{\lambda T}}{\lambda^{p}} \int_0^T \frac{1}{x^{u-1}} dx \quad ( \int_0^T \frac{1}{x^{u}} dx)^{p-1} \nonumber \\
= & \frac{\Gamma(u)^p e^{\lambda T}}{\lambda^{p}} \frac{T^{(1-u)p+1}}{(1-u)^{p-1}(2-u)},
\label{eq:HH}
\end{align}
so that
\begin{align}
\label{eq:B1}
B_1(\tilde{k}^1,..,\tilde{k}^P)  \le  \frac{\Gamma(u)^p e^{\lambda T}}{\lambda^{p}} \frac{T^{(1-u)p+1+2\theta}}{(1-u)^{p-1}(2-u)}C(\sigma)^{p-1} \hat K^2  K^{2p}. 
\end{align}
In a similar way introducing $\tilde{C}(\sigma) =\E(\phi_k^2)^{\frac{1}{2}}$ and using the tower rule, the independence of the different random variables  and Cauchy Schwartz,
\begin{align}
\label{eq:B2}
B_2(\tilde{k}^1,..,\tilde{k}^{i+1})  &= 
\E \big( 1_{T_{\tilde{k}^{i+1}} < T}
\E( \hat f_{\tilde{k}^{i+1}}^2 \prod_{j=2}^{i+1} \phi_{\tilde k^j} / \tau_{\tilde{k}^1},..,\tau_{\tilde{k}^{i+1}})
\frac{1}{\rho(\tau_{\tilde{k}^1})^2}
\prod_{j=2}^{i+1} \frac{1}{\rho(\tau_{\tilde{k}^j})^2 \tau_{\tilde{k}^j}  }
\big)  \nonumber\\
& \le \sup_{t \in [0,T]} \E[ f(t,X_t,Du(t,X_t))^4]^{\frac{1}{2}} \E(\prod_{j=2}^{i+1} \phi_{\tilde k^j}^2)^{\frac{1}{2}} E \big( 1_{T_{\tilde{k}^{i+1}} < T} \frac{1}{\rho(\tau_{\tilde{k}^1})^2}\prod_{j=2}^{i+1} \frac{1}{\rho(\tau_{\tilde{k}^j})^2 \tau_{\tilde{k}^j} }
\big)  \nonumber\\
& \le 
\frac{\Gamma(u)^{i+1}e^{\lambda T}}{\lambda^{i+1}} \frac{T^{(1-u)(i+1)+1}}{(1-u)^{i}(2-u)} \tilde C(\sigma)^{i} \sup_{t \in [0,T]} \E[ f(t,X_t,Du(t,X_t))^4]^{\frac{1}{2}}.
\end{align}
For $B_3$ we divide the calculation into two cases :
\begin{align}
B_3(\tilde{k}^1) = & \E( 1_{\tau_{\tilde{k}^1}>T} \frac{1}{\bar F(T)^2}  \E(g(X_T)^2)) \nonumber \\
= & \frac{1}{\bar F(T)}  \E(g(X_T)^2) \nonumber \\
= & \frac{\Gamma(u)}{\Gamma(u)-\gamma(u,\lambda T)}  \E(g(X_T)^2).
\label{eq:B31}
\end{align}
For $i>0$, using the Lipschitz property of $g$, the tower property, the independence of the different random variables, and noting $$\bar C(\mu, \sigma, T)= \E[ (2 ||\mu||_2^2 T + 2 || \sigma G_{\tilde{k}^{i+1}}||_2^2) \phi_{\tilde{k}^{i+1}} ],$$ we have
the following bound:
\begin{align}
B_3(\tilde{k}^1, ..., \tilde{k}^{i+1}) \le & \E \big [1_{T_{\tilde{k}^{i+1} > T}} 1_{T_{\tilde{k}^{i} < T}} \tilde K^2  \frac{(2 ||\mu||_2^2 T + 2 || \sigma G_{\tilde{k}^{i+1}}||_2^2) \phi_{\tilde{k}^{i+1}}}{ \bar F(T)^2} \frac{C(\sigma)^{i-1}}{\rho(\tau_{\tilde{k}^1})^2} \prod_{j=2}^{i}  \frac{1}{\rho(\tau_{\tilde{k}^j})^2 \tau_{\tilde{k}^j}} \big ]\nonumber \\
\le & \frac{ \bar C(\mu, \sigma,T) \tilde K^2}{\bar F(T)^2} C(\sigma)^{i-1} \frac{\Gamma(u)^i e^{\lambda T}}{\lambda^{i}} \int_0^T \frac{1}{x^{u-1}} dx \quad ( \int_0^T \frac{1}{x^{u}} dx)^{i-1} \int_T^\infty \rho(x) dx \nonumber \\
= &  C(\sigma)^{i-1} \frac{\Gamma(u)^{i+1}}{\lambda^{i}} \frac{T^{(1-u)i+1}}{(1-u)^{i-1}(2-u)} \frac{ \bar C(\mu, \sigma,T) \tilde K^2 e^{\lambda T}}{\Gamma(u)-\gamma(u,\lambda T )}. 
\label{eq:B32}
\end{align}
Plugging \eqref{eq:B1}, \eqref{eq:B2}, \eqref{eq:B31}, \eqref{eq:B32} into \eqref{eq:A1}, \eqref{eq:A2} and \eqref{eq:A3} that we insert into \eqref{est::finalErroSemi} gives the desired estimation.
\qquad \end{proof}

\subsection{A second estimator}
In this section we present a scheme derived from \cite{warin2017variations}.\\
We extend the notations from the previous section  \ref{sec:fistEstim}.
Let set $p \in \N^{*}$.
We construct the sets $Q^o_i$  for $i = 1, .., p$, such that
$Q^o_1= \{ (k_1), (k_1^{-}) \}$ where $k_1 \in\{1,..,N_1\}$, so that to a particle noted $(k_1)  \in Q_1$, we associate a antithetic particle noted $k_1^{-}$. 
Then the set $Q^o_i$ are defined by recurrence : 
\begin{flalign*}
Q^o_{i+1} = \{ (k_1,..,k_i, k_{i+1}) / (k_1,..,k_i) \in Q^o_{i}, k_{i+1} \in \{ 1,..,N_{i+1}, 1^{-},.., N^{-}_{i+1} \} \}
\end{flalign*}
To a particle $k=(k_1, .., k_i) \in Q^o_i$ we associate its original particle $o(k) \in Q_i$ such that $o(k)= (\hat k_1,.. \hat k_i)$ where $\hat k_j = l$ if $k_j = l$ or $l^{-}$.
Further, when $k = (k_1, \cdots, k_i)$ is such that $k_i \in \N$, we denote $k^{-} := (k_1, \cdots, k_{i-1}, k_i^{-})$.\\
By convention $T_k = T_{o(k)}$, $ \tau_k = \tau_{o(k)}$ and $\bar W_t^k =\bar W_t^{o(k)}$.
For $k=(k_1, ...,k_i) \in Q^o_i$ we introduce the set 
\begin{itemize}
\item $\tilde Q^o(k) = \{ l =(k_1,..,k_i, m )/  m \in \{1,.., N_i\} \} \subset Q_{i+1}^o$
\item and
$\hat Q^o(k) = \{ l =(k_1,..,k_i, m )/  m \in \{1,.., N_i , 1^{-},.., N_i^{-}\} \} \subset Q_{i+1}^o$
\end{itemize}
For $k = (k_1,..,k_i) \in Q^o_i$  and $\tilde k =(k_1,..,k_i, k_{i+1})  \in \hat Q^o(k)$ we define the following trajectories :
	\begin{flalign}\label{eq:brownRenorm}
		W^{\tilde k}_s
		~:=~&
		W^{k}_{T_{k}}
		~+~
		\1_{k_{i+1} \in \N}
		\bar W^{o(\tilde k)}_{s - T_{k}} 
        ~-~
		\1_{k_{i+1} \notin \N}
		\bar W^{o(\tilde k)}_{s - T_{k}},
		~~~\mbox{and}~~\\
		X^{\tilde k}_s := & x+ \mu s +\sigma  W^{\tilde k}_{s},
		~~~\forall s \in [T_{k}, T_{\tilde k}] .
	\end{flalign}

%In order to understand what these different trajectories represent, suppose that  $d=1$, $\mu= 0$, $\sigma =1$ and let us consider the original particle $k= (1,1,1)$  such that $T_{(1,1,1)}=T$.\\
%Following equation \eqref{eq:brownRenorm},
%\begin{align*}
%X^{(1,1,1)}_T = & \hat W^{(1)}_{T_{(1)}} + \hat W^{(1,1)}_{T_{(1,1)}-T_{(1)}} + \hat W^{(1,1,1)}_{T-T_{(1,1)}}, \\
%X^{(1^{-},1,1)}_T= & -\hat W^{(1)}_{T_{(1)}} + \hat W^{(1,1)}_{T_{(1,1)}-T_{(1)}} + \hat W^{(1,1,1)}_{T-T_{(1,1)}}, \\
%X^{(1,1^{-},1)} = & \hat W^{(1)}_{T_{(1)}} - \hat W^{(1,1)}_{T_{(1,1)}-T_{(1)}} + \hat W^{(1,1,1)}_{T-T_{(1,1)}}, \\
%X^{(1^{-},1^{-},1)}_T = & -\hat W^{(1)}_{T_{(1)}} - \hat W^{(1,1)}_{T_{(1,1)}-T_{(1)}} + \hat W^{(1,1,1)}_{T-T_{(1,1)}}, \\
%X^{(1^{-},1^{-},1^{-})}_T = & -\hat W^{(1)}_{T_{(1)}} - \hat W^{(1,1)}_{T_{(1,1)}-T_{(1)}} - \hat W^{(1,1,1)}_{T-T_{(1,1)}}, \\
%X^{(1^{-},1,1^{-})}_T = & -\hat W^{(1)}_{T_{(1)}} + \hat W^{(1,1)}_{T_{(1,1)}-T_{(1)}} - \hat W^{(1,1,1)}_{T-T_{(1,1)}}, \\
%...&
%\end{align*}
%Then for example $X^{(1,1,1)}_T$ and $X^{(1,1,1^{-})}_T$ only differ from $2 \hat W^{(1,1,1)}_{T_{(1)}}$ so that if $T_{(1,1,1)} - T_{1,1}$ is small, the trajectories $ X^{(1,1,1)}_t$
%and $X^{(1,1,1^{-})}_t$ are very close, leading to values of  $D\ bar u ^p_{(1,1,1)}-D\ bar u ^p_{(1,1,1^{-})}$ in $o(\sqrt{ T_{(1,1,1)} - T_{1,1}})$.
%\\
Using the previous definitions,
we consider the estimator defined by:
\begin{align}
\left \{
 \begin{array}{ll}
\bar u_0^p = & \frac{1}{N_0} \sum_{j=1}^{N_0}  \frac{\big(  \hat \phi\big( 0, T_{(j)}, X^{(j)}_{T_{(j)}}, D \bar u_{(j)}^p\big) +
\hat \phi\big( 0,T_{(j)}, X^{(j)^{-}}_{T_{(j)}}, D \bar u_{(j)^{-}}^p\big)\big)}{2}
\big),\\
D\bar u_{k}^p = &  \frac{1}{N_i} \sum_{\tilde k \in \tilde Q^o(k)}  \tilde W^{\tilde k} \frac{1}{2} \big( \hat \phi\big(T_k,T_{\tilde k},X^{\tilde k}_{T_{\tilde k}}, D\bar u_{\tilde k}^p\big) -\hat \phi\big(T_k,T_{\tilde k},  X^{\tilde k^{-}}_{T_{\tilde k}}, D\bar u_{\tilde k^{-}}^p\big) \big), \\
& \quad \mbox{ for }  k =( k_1,..k_i) \in Q^o_i, i <p,  \\
D\bar u_{\tilde k}^p = & D g(X^{\tilde k}_{T_{\tilde k}}) \quad  \mbox{for } \tilde k \in Q_p^o
\end{array}
\right.
\label{eq:estimSemi2}
\end{align}
 where $\hat \phi $ is defined by equation \eqref{eq:phihat} and $ \tilde W^{\tilde k} = \sigma^{-\top} \frac{\bar W_{ T_{\tilde k}- T_k}^{o(\tilde k)}}{T_{\tilde k}- T_k}$ .\\
The idea is that, for a given $k =(k_1,..,k_i) \in Q^o_i$ and a given $\tilde k \in \tilde Q^o(k)$, if we have
$\bar W_{ T_{\tilde k}- T_k}$ very small then
\begin{align*}
\frac{1}{2} \big( \hat \phi\big(T_k,T_{\tilde k},X^{\tilde k}_{T_{\tilde k}}, D\bar u_{\tilde k}^p\big) -  \hat \phi\big(T_k,T_{\tilde k},  X^{\tilde k^{-}}_{T_{\tilde k}}, D\bar u_{\tilde k^{-}}^p)\big) \simeq & D_u\hat \phi\big(T_k,T_{\tilde k},X^{\tilde k}_{T_{\tilde k}}, D\bar u_{\tilde k}^p\big) (D\bar u_{\tilde k}^p - D\bar u_{\tilde k^{-}}^p\big) + \\
&  D_x\hat \phi\big(T_k,T_{\tilde k},X^{\tilde k}_{T_{\tilde k}}, D\bar u_{\tilde k}^p\big) (X^{\tilde k}_{T_{\tilde k}} -X^{\tilde k^{-}}_{T_{\tilde k}}) \\
 \simeq & C \sqrt{\tau^{\tilde k}}.
\end{align*}
So  $\tilde W^{\tilde k}  \frac{1}{2} \big( \hat \phi\big(T_k,T_{\tilde k},X^{\tilde k}_{T_{\tilde k}}, D\bar u_{\tilde k}^p\big) -  \hat \phi\big(T_k,T_{\tilde k},  X^{\tilde k^{-}}_{T_{\tilde k}}, D\bar u_{\tilde k^{-}}^p)\big)$, as function of  $\tau^{\tilde k}$,  should be bounded, and  we hope that estimator \eqref{eq:estimSemi2} has a much smaller variance than estimator \eqref{eq:estimSemi1}.
\begin{remark}
As in a previous scheme an estimation of the gradient  is obtained as:
\begin{flalign*}
D \bar u_0^p = & \frac{1}{N_0} \sum_{j=1}^{N_0}  \frac{1}{2} \sigma^{-\top} \frac{W_{T_{(j)}}}{T_{(j)}}
\big( \hat \phi\big( 0,T_{(j)}, X^{(j)}_{T_{( j)}}, D\bar u_{(j)}^p\big) -
 \hat \phi\big( 0,T_{(j)}, X^{(j)^{-}}_{T_{(j)}}, D\bar u_{(j)^{-}}^p\big) 
\big)
\end{flalign*} 
\end{remark}
We need other  assumptions on the solution  and the driver to fully exploit this scheme:
 \begin{assumption}
 \label{ass::DuLipSemi}
 $Du$ is uniformly Lipschitz in $x$  such that for $\bar K >0$:
 \begin{align*}
 || Du(t,x) -Du(t,y)||_2 \le \bar K ||x-y||_2 \quad \forall (t,x,y) \in [0,T]\times \R^d \times \R^d.
 \end{align*}
 \end{assumption}
 \begin{assumption}
 \label{ass::fxLipSemi}
$f$  is uniformly Lipschitz in $x$ such that there exists $\underbar K > 0$
 \begin{align*}
 |f(t,x,z) -f(t,y,z)| \le \underbar{K} ||x-y||_2, \quad  \forall (t,x,y,z) \in [0,T]\times \R^d \times \R^d \times \R^d.
 \end{align*}
 \end{assumption}
We now give the error estimate with the second scheme
\begin{proposition}
\label{propGhost}
Suppose that $\rho$ is the density of a gamma law so that $ \rho(x)=\lambda^u x^{u-1} \frac{e^{-\lambda x}}{\Gamma(u)}$,  suppose that  $u<1$, and suppose that  assumptions \ref{ass::lipFSemi},  \ref{ass::uHolderSemi},  \ref{ass::gLipSemi},\ref{ass::DuLipSemi} and \ref{ass::fxLipSemi} are satisfied then  there are two  functions $C$ and $\tilde C$ depending on $\sigma$ and one $\bar C$ depending on $K$, $\bar K$, $\underbar K$ and  $\sigma$  such that using estimator \eqref{eq:estimSemi2} to solve equation  \eqref{eqPDESemi}, we  have the error estimate
\begin{align}
\label{est::final}
 \E\big((\bar u_{\emptyset}^p - u(0,x))^2 \big) \le
 \prod_{i=1}^p (1+\frac{8}{N_{i-1}}) \frac{\Gamma(u)^p e^{\lambda T}}{\lambda^{p}} \frac{T^{(1-u)p+1+2\theta}}{(1-u)^{p-1}(2-u)}C(\sigma)^{p-1} \hat K^2  K^{2p} + \nonumber \\
   \sum_{i=1}^{p-1} \frac{K^{2i}}{N_{i}} \prod_{j=1}^i (1+\frac{8}{N_{j-1}}) \bar C( \sigma,K,\bar K , \underbar K ) C(\sigma)^{i-1}
 \frac{\Gamma(u)^{i+1}e^{\lambda T}}{\lambda^{i+1}} \frac{T^{(1-u)i+3-u}}{(2-u)^2(1-u)^{i-1}}  +  \nonumber  \\
 \sum_{i=1}^{p-1} \frac{K^{2i}}{2 N_{i}} \prod_{j=1}^i (1+\frac{8}{N_{j-1}})
 C( \sigma)^{i-1} \frac{\tilde{K}^2 \Gamma(u)^{i+1} e^{\lambda T}}{ \lambda^i  (\Gamma(u)-\gamma(u,\lambda T ))} \bar C(\sigma) \frac{T^{(1-u)i+1}}{(2-u)(1-u)^{i-1}}
 + \nonumber \\
 \frac{4}{N_0} \frac{\Gamma(u)}{\lambda} e^{\lambda T} \frac{T^{2-u}}{2-u} \hat F +\frac{2}{ N_0} \frac{\Gamma(u)}{\Gamma(u)-\gamma(u,\lambda T)}   \E(g(X_T)^2)
\end{align}
with $\hat F$ given by equation \eqref{eq:FHat}.
\end{proposition}

\begin{proof}
First notice that under assumption \ref{ass::uHolderSemi},   $u$ satisfies \eqref{eq:recSemi} and 
then, for $j \in [0,N_0]$
\begin{flalign*}
 u(0,x) = 
 \E_{0,x} \big[ \frac{ \hat \phi\big( 0, T_{(j)}, X^{(j)}_{T_{(j)}}, Du(T_{(j)},X^{(j)}_{T_{(j)}})\big) +
\hat \phi\big( 0,T_{(j)}, X^{(j)^{-}}_{T_{(j)}}, Du(T_{(j)},X^{(j)^{-}}_{T_{(j)}}\big)}{2}\big],
\end{flalign*}
and because $u$ satisfies  equation \eqref{eq:ruelDu}, we have that for $k \in Q^o_i$, $\tilde k \in \tilde Q^o(k)$,
\begin{flalign*}
D u(T_k,X^k_{t_k})= &  \E_{T_k,X^k_{t_k}} \big(  \sigma^{-\top} \frac{\bar W_{ T_{\tilde k}- T_k}^{o(\tilde k)}}{T_{\tilde k}- T_k}\frac{1}{2} \big( \hat \phi\big(T_k,T_{\tilde k},X^{\tilde k}_{T_{\tilde k}}, D u(T_{\tilde k},X^{\tilde k}_{T_{\tilde k}})\big) -  \hat \phi\big(T_k,T_{\tilde k},  X^{\tilde k^{-}}_{T_{\tilde k}}, D u(T_{\tilde k},  X^{\tilde k^{-}}_{T_{\tilde k}})\big) \big)
\end{flalign*}

Then we can introduce  for $ k \in Q_i$, $1 \le i < p$:
\begin{flalign*}
E_k := & \E_{T_{k},X^{k}_{T_{k}} }\big( || D \bar u_{k}^p - Du(T_k,X^k_{T_k})||^2_2  1_{T_k <T} \big).
\end{flalign*} 
In exactly the same way as in the demonstration of proposition \ref{theo2}, we have the following result similar to the one given by equation \eqref{eq:ErrLessD}:
\begin{align}
E_k \le& \frac{1}{N_i} (1+\frac{8}{ N_i}) \sum_{\tilde k \in \tilde Q^o(k)} E_{T_k,X^{k}_{T_{k}}} \big(  1_{T_{\tilde k < T}} \frac{(\bar W^{\tilde{k}}_{\tau_{\tilde{k}}})^\top \sigma^{-1} \sigma^{-\top} \bar W^{\tilde{k}}_{\tau_{\tilde{k}}}}{\tau_{\tilde{k}}^2} \frac{1}{4\rho(\tau_{\tilde k})^2} \big( f(T_{\tilde k},X^{\tilde k}_{T_{\tilde k}},D \bar u_{\tilde{k}}^p)- \nonumber \\ 
& f(T_{\tilde k},X^{\tilde k^{-}}_{T_{\tilde k}},D \bar u_{\tilde{k}^{-}}^p) - f(T_{\tilde k},X^{\tilde k}_{T_{\tilde k}},Du(T_{\tilde k},X^{\tilde{k}}_{T_{\tilde{k}}})) +
f(T_{\tilde k},X^{\tilde k^{-}}_{T_{\tilde k}},Du(T_{\tilde k},X^{\tilde{k}^{-}}_{T_{\tilde{k}}})))^2\big) + \nonumber
\\
& 4 \frac{1}{N_i^2} \sum_{\tilde k \in \tilde Q^o(k)} \E_{T_k,X^k_{T_k} }\big(  1_{T_{\tilde k < T}} \frac{(\bar W^ {\tilde{k}}_{\tau_{\tilde{k}}})^\top \sigma^{-1} \sigma^{-\top} \bar W^{\tilde{k}}_{\tau_{\tilde{k}}}}{\tau_{\tilde{k}}^2} \nonumber \\
& (\frac{f(T_{\tilde k}, X^{\tilde k}_{T_{\tilde k}}, Du(T_{\tilde k},X^{\tilde k}_{T_{\tilde k}}))-f(T_{\tilde k}, X^{\tilde k^{-}}_{T_{\tilde k}}, Du(T_{\tilde k},X^{\tilde k^{-}}_{T_{\tilde k}}))}{2\rho(\tau_{\tilde k})} )^2\big) +  \nonumber \\
&2 \frac{1}{N_i^2} \sum_{\tilde k \in \tilde Q^o(k)} \E_{T_k,X^k_{T_k} }\big( 1_{T_{\tilde k\ge T}} \frac{(\bar W^{\tilde{k}}_{T-T_k})^\top \sigma^{-1} \sigma^{-\top} \bar W^{\tilde{k}}_{T-T_k}}{(T-T_k)^2} \frac{(g(X^{\tilde k}_T) -g(X^{\tilde k^{-}}_T))^2}{ 4 \Fb(T-T_k)^2} \big).
\end{align}
Using $ (a+b)^2 \le 2a^2 +2b^2$, the Lipschitz property of $f$, we get

\begin{align}
E_k \le& \frac{1}{2N_i} (1+\frac{8}{ N_i}) \sum_{\tilde k \in \hat Q^o(k)} E_{T_k,X^{k}_{T_{k}}} \big(  \frac{(\bar W^{\tilde{k}}_{\tau_{\tilde{k}}})^\top \sigma^{-1} \sigma^{-\top} \bar W^{\tilde{k}}_{\tau_{\tilde{k}}}}{\tau_{\tilde{k}}^2} \frac{K^2}{\rho(\tau_{\tilde k})^2} E_{\tilde k} \big) + \nonumber
\\
&  \frac{1}{ N_i^2} \sum_{\tilde k \in \tilde Q^o(k)} \E_{T_k,X^k_{T_k} }\big(  1_{T_{\tilde k < T}} \frac{(\bar W^{\tilde{k}}_{\tau_{\tilde{k}}})^\top \sigma^{-1} \sigma^{-\top} \bar W^{\tilde{k}}_{\tau_{\tilde{k}}}}{\tau_{\tilde{k}}^2} \nonumber \\
& (\frac{f(T_{\tilde k}, X^{\tilde k}_{T_{\tilde k}}, Du(T_{\tilde k},X^{\tilde k}_{T_{\tilde k}}))-f(T_{\tilde k}, X^{\tilde k^{-}}_{T_{\tilde k}}, Du(T_{\tilde k},X^{\tilde k^{-}}_{T_{\tilde k}}))}{\rho(\tau_{\tilde k})} )^2\big) +  \nonumber \\
& \frac{1}{2 N_i^2} \sum_{\tilde k \in \tilde Q^o(k)} \E_{T_k,X^k_{T_k} }\big( 1_{T_{\tilde k\ge T}} \frac{(\bar W^{\tilde{k}}_{T-T_k})^\top \sigma^{-1} \sigma^{-\top} \bar W^{\tilde{k}}_{T-T_k}}{(T-T_k)^2} \frac{(g(X^{\tilde k}_T) -g(X^{\tilde k^{-}}_T))^2}{ \Fb(T-T_k)^2} \big).
\label{eq:ErrLessD2}
\end{align}

Similarly to \eqref{eq:E0}
\begin{align}
 D =& \E\big((\bar u_{\emptyset}^p - u(0,x))^2 \big)  \nonumber \\
 \le& \frac{1}{N_0} (1+\frac{8}{N_0}) \sum_{\tilde k \in \tilde Q^o(\emptyset)} \E \big( 1_{T_{\tilde k < T}}  \frac{K^2}{4 \rho(\tau_{\tilde k})^2}  \big( f(T_{\tilde k},X^{\tilde k}_{T_{\tilde k}},D \bar u_{\tilde{k}}^p)+\\ \nonumber
& f(T_{\tilde k},X^{\tilde k^{-}}_{T_{\tilde k}},D \bar u_{\tilde{k}^{-}}^p) -  f(T_{\tilde k},X^{\tilde k}_{T_{\tilde k}},Du(T_{\tilde k},X^{\tilde{k}}_{T_{\tilde{k}}})) -
f(T_{\tilde k},X^{\tilde k^{-}}_{T_{\tilde k}},Du(T_{\tilde k},X^{\tilde{k}^{-}}_{T_{\tilde{k}}})))^2\big)\big) + \nonumber
\\
& 4 \frac{1}{N_0^2} \sum_{\tilde k \in \tilde Q^o(\emptyset)} \E \big(  1_{T_{\tilde k < T}} (\frac{f(T_{\tilde k},X^{\tilde k}_{T_{\tilde k}},  Du(T_{\tilde k},X^{\tilde k}_{T_{\tilde k}}))+f(T_{\tilde k},X^{\tilde k^{-}}_{T_{\tilde k}},Du(T_{\tilde k},X^{\tilde{k}^{-}}_{T_{\tilde{k}}})) }{2\rho(\tau_{\tilde k})} )^2\big) +  \nonumber \\
&2 \frac{1}{N_0^2} \sum_{\tilde k \in \tilde Q^o(\emptyset)} \E_{T_k,X^k_{T_k} }\big( 1_{T_{\tilde k\ge T}} \frac{(g(X^{\tilde k}_T)+g(X^{\tilde k^{-}}_T))^2}{ 4\Fb(T-T_k)^2} \big)
\end{align}
so that using the Lipschitz property of $f$
\begin{align}
 \E\big((\bar u_{\emptyset}^p - u(0,x))^2 \big) \le& \frac{1}{2N_0} (1+\frac{8}{N_0}) \sum_{\tilde k \in \hat Q^o(\emptyset)} \E \big( 1_{T_{\tilde k < T}}  \frac{K^2}{ \rho(\tau_{\tilde k})^2}   E_{\tilde{k}} \big) + \nonumber
\\
& 2 \frac{1}{N_0^2} \sum_{\tilde k \in \hat Q^o(\emptyset)} \E \big(  1_{T_{\tilde k < T}} (\frac{f(T_{\tilde k},X^{\tilde k}_{T_{\tilde k}},  Du(T_{\tilde k},X^{\tilde k}_{T_{\tilde k}}))}{\rho(\tau_{\tilde k})} )^2\big) +  \nonumber \\
& \frac{1}{ N_0^2} \sum_{\tilde k \in \hat Q^o(\emptyset)} \E_{T_k,X^k_{T_k} }\big( 1_{T_{\tilde k\ge T}} \frac{g(X^{\tilde k}_T)^2}{\Fb(T-T_k)^2} \big).
\end{align}
Then we get that
\begin{align}
\label{est::finalH}
 \E\big((\bar u_{\emptyset}^p - u(0,x))^2 \big) \le&  \hat A_1 + \hat A_2 +  \hat A_3,
\end{align}
where the terms $\hat A_1$ , $\hat A_2$ and $\hat A_3$ are given by:
\begin{itemize}
  \item
\begin{align}
\hat A_1 = \prod_{i=1}^p \frac{1}{2 N_{i-1}^2} (1+\frac{8}{ N_{i-1}}) \sum_{\tilde k^1 \in \hat Q^o(\emptyset)} ... \sum_{\tilde k^p \in \hat Q^o(\tilde{k}^{p-1})}  B_1(\tilde{k}^1,..,\tilde{k}^P)  
\label{eq:HA1}
\end{align}
where $B_1$ is given by \eqref{eq:alllB} and bounded by \eqref{eq:B1},
\item
\begin{flalign}
\hat A_2 = & \frac{2}{N_0^2} \sum_{\tilde k^1 \in \hat Q^o(\emptyset)}  B_2(\tilde{k}^1) +  \nonumber \\
&\sum_{i=1}^{p-1} \frac{K^{2i}}{ N_{i}^2} \prod_{j=1}^i \frac{1}{2 N_{j-1}} (1+\frac{8}{N_{j-1}})
\sum_{\tilde k^1 \in \hat Q^o(\emptyset)} ...\sum_{\tilde k^{i} \in \hat Q^o(\tilde{k}^{i-1})}  \sum_{\tilde k^{i+1} \in \tilde Q^o(\tilde{k}^{i})}  \hat B_2(\tilde{k}^1,..,\tilde{k}^{i+1}),
\label{eq:HA2}
\end{flalign}
where $B_2$ is given by equation \eqref{eq:alllB} and bounded by equation \eqref{eq:B2} so that
\begin{align}
  \label{eq:B2One}
  B_2(\tilde{k}^1) \le \frac{\Gamma(u)}{\lambda} e^{\lambda T} \frac{T^{2-u}}{2-u} \hat F,
\end{align}
with $\hat F$ given by equation \eqref{eq:FHat},
and where using notation given by \eqref{eq:not1}
\begin{align}
  \label{eq:HatB1}
\hat B_2(\tilde{k}^1,..,\tilde{k}^{i+1}) = & \E \big [ 1_{T_{\tilde{k}^{i+1} < T}} (\hat f_{\tilde k^{i+1}}-\hat f_{(\tilde k^{i+1})^{-}})^2
  \prod_{j=1}^{i+1} \frac{1}{\rho(\tau_{\tilde{k}^j})^2} \prod_{j=2}^{i+1} \psi_{\tilde{k}^j} \Big],
\end{align}
\item 
\begin{align}
\hat A_3 =& \frac{1}{N_0^2} \sum_{\tilde k^1 \in \hat Q^o(\emptyset)} \hat B_{3}(\tilde{k}^1) + \nonumber \\
 &\sum_{i=1}^{p-1} \frac{K^{2i}}{ N_{i}^2} \prod_{j=1}^i \frac{1}{2 N_{j-1}} (1+\frac{8}{N_{j-1}})
\sum_{\tilde k^1 \in \hat Q^o(\emptyset)} ... \sum_{\tilde k^{i} \in \hat Q^o(\tilde{k}^{i-1})}  \sum_{\tilde k^{i+1} \in \tilde Q^o(\tilde{k}^{i})}  \hat B_{3}(\tilde{k}^1,..,\tilde{k}^{i+1}),
\label{eq:HA3}
\end{align}
where using notation given by \eqref{eq:not1}
\begin{flalign}
  \label{eq:HatB2}
 \hat B_{3}(\tilde{k}^1,..,\tilde{k}^{i+1}) = &\E [1_{T_{\tilde{k}^{i+1} > T}} 1_{T_{\tilde{k}^{i} < T}} \frac{(g(X^{\tilde k ^{i+1}}_{T})- 1_{i>1} g(X^{(\tilde k ^{i+1})^{-}}_{T_{\tilde{k}^i}}))^2}{ \bar  F(T-T_{\tilde{k}^{i}})^2} 
 \prod_{j=1}^{i} \frac{1}{\rho(\tau_{\tilde{k}^j})^2} \prod_{j=2}^{i} \psi_{\tilde{k}^j}  \nonumber \\
 & \frac{(\bar W^{\tilde k^{i+1}}_{T-T_{\tilde k^i}})^\top \sigma^{-1} \sigma^{-\top} \bar W^{\tilde k^{i+1} }_{T-T_{\tilde k^i}}}
 {(T-T_{\tilde k^i})^2} ]
 %\prod_{j=1}^{i} \frac{1}{\rho(\tau_{\tilde{k}^j})^2} \prod_{j=2}^{i+1} \psi_{\tilde{k}^j}\big ].
 %\label{eq:alllB}
\end{flalign}
\end{itemize}
We can bound $\hat B_2$ using definitions in equations \eqref{eq:not1} and\eqref{eq::phi}
\begin{flalign*}
\hat B_2(\tilde{k}^1,..,\tilde{k}^{i+1}) = & \E \big [ 1_{T_{\tilde{k}^{i+1} < T}}  \E\big( (\hat f_{\tilde k^{i+1}}-\hat f_{(\tilde k^{i+1})^{-}})^2  \prod_{j=2}^{i+1} \phi_{\tilde{k}^{j}} | \tau_{\tilde{k}^{1}},...,\tau_{\tilde{k}^{i+1}}\big) \\
&
 \frac{1}{\tau_{\tilde{k}^{i+1}} \rho(\tau_{\tilde{k}^{i+1}})^2}\prod_{j=1}^{i} \frac{1}{\rho(\tau_{\tilde{k}^j})^2} \prod_{j=2}^{i} \frac{1}{\tau_{\tilde{k}^j}} \Big] .
\end{flalign*}
Using  assumptions \ref{ass::lipFSemi}, \ref{ass::DuLipSemi} and \ref{ass::fxLipSemi},
\begin{align*}
|\hat f_{\tilde k^{i+1}}-\hat f_{(\tilde k^{i+1})^{-}}| \le & 2 \bar  K || \sigma W^{\tilde{k}^{i+1}}_{T_{\tilde{K}^{i+1}}-T_{\tilde k^i}}||_2
+  K ||Du(T_{\tilde{K}^{i+1}}, X^{\tilde{k}^{i+1}}_{T_{\tilde{K}^{i+1}}})-
Du(T_{\tilde{K}^{i+1}}, X^{(\tilde{k}^{i+1})^{-}}_{T_{\tilde{K}^{i+1}}})||_2 \\
\le & 2 ( \underbar  K + K \bar K ) || \sigma W^{\tilde{k}^{i+1}}_{T_{\tilde{K}^{i+1}}-T_{\tilde k^i}}||_2 
\end{align*}
so that  :
\begin{align*}
R =&   \E\big( (\hat f_{\tilde k^{i+1}}-\hat f_{(\tilde k^{i+1})^{-}})^2  \prod_{j=2}^{i+1} \phi_{\tilde{k}^{j}} |  \tau_{\tilde{k}^{1}},...,\tau_{\tilde{k}^{i+1}}\big) \\
& \le \bar C( \sigma,K,\bar K , \underbar K ) C(\sigma)^{i-1} \tau_{\tilde{k}^{i+1}}.
\end{align*}
Then
\begin{flalign}
\hat B_2(\tilde{k}^1,..,\tilde{k}^{i+1}) \le & \bar C( \sigma,K,\bar K , \underbar K ) C(\sigma)^{i-1} \E \big [ 1_{T_{\tilde{k}^{i+1} < T}}  
 \frac{1}{ \rho(\tau_{\tilde{k}^{i+1}})^2}\prod_{j=1}^{i} \frac{1}{\rho(\tau_{\tilde{k}^j})^2} \prod_{j=2}^{i} \frac{1}{\tau_{\tilde{k}^j}} \Big] \nonumber  \\
 \le & \bar C( \sigma,K,\bar K , \underbar K ) C(\sigma)^{i-1}
 \frac{\Gamma(u)^{i+1} e^{\lambda T}}{\lambda^{i+1}} (\int_0^T\frac{1}{x^{u-1}} dx)^2
 (\int_0^T\frac{1}{x^{u}}dx)^{i-1}, \nonumber \\
 = & \bar C( \sigma,K,\bar K , \underbar K ) C(\sigma)^{i-1}
 \frac{\Gamma(u)^{i+1}e^{\lambda T}}{\lambda^{i+1}} \frac{T^{(1-u)i+3-u}}{(2-u)^2(1-u)^{i-1}}.
 \label{eq:HB2}
\end{flalign}
Similarly
\begin{align}
\hat B_3(\tilde{k}^1) = B_3(\tilde{k}^1)
\label{eq:HB31}
\end{align}
with $B_3$ given by \eqref{eq:B31}.\\
For the general $B_3$ term for $i>0$, using assumption \ref{ass::gLipSemi} and  using the notation $ \bar C(\sigma)=\E[4\phi_{\tilde{k}^{i+1}} || \sigma G_{\tilde k^{i+1}} ||_2^2]$:
\begin{flalign}
\hat B_{3}(\tilde{k}^1,..,\tilde{k}^{i+1}) = &\E \big [1_{T_{\tilde{k}^{i+1} > T}} 1_{T_{\tilde{k}^{i} < T}} \frac{(g(X^{\tilde k ^{i+1}}_{T})-g(X^{(\tilde k ^{i+1})^{-}}_{T_{\tilde{k}^i}}))^2}{  (T-T_{\tilde k^{i}}) \bar  F(T-T_{\tilde{k}^{i}})^2}   \frac{1}{\rho(\tau_{\tilde{k}^1})^2} \prod_{j=2}^{i} \frac{1}{   \tau_{\tilde{k}^j}   \rho(\tau_{\tilde{k}^j})^2} \prod_{j=2}^{i+1} \phi_{\tilde{k}^j}\big ] \nonumber \\
\le & C( \sigma)^{i-1} \frac{\tilde{K}^2}{\bar F(T)^2} \E[4\phi_{\tilde{k}^{i+1}} || \sigma G_{\tilde k^{i+1}} ||_2^2]  \E \big (1_{T_{\tilde{k}^{i+1} > T}} 1_{T_{\tilde{k}^{i} < T}} \frac{1}{\rho(\tau_{\tilde{k}^1})^2} \prod_{j=2}^{i}   \frac{1}{   \tau_{\tilde{k}^j}      \rho(\tau_{\tilde{k}^j})^2}  \big) \nonumber \\
\le & C( \sigma)^{i-1} \frac{\tilde{K}^2 e^{\lambda T} \Gamma(u)^i }{ \lambda^i \bar F(T)^2} \bar C(\sigma) \int_0^T \frac{1}{x^{u-1}}dx \big( \int_0^T \frac{1}{x^u} dx \big)^{i-1} \int_T^\infty \rho(x) dx \nonumber\\
\le & C( \sigma)^{i-1} \frac{\tilde{K}^2 \Gamma(u)^i e^{\lambda T}}{\bar F(T) \lambda^i} \bar C(\sigma) \frac{T^{(1-u)i+1}}{(2-u)(1-u)^{i-1}}
\label{eq:HB32}
\end{flalign}

Plugging \eqref{eq:B1}, \eqref{eq:B2One}, \eqref{eq:HB2},\eqref{eq:HB31},\eqref{eq:HB32} in \eqref{eq:HA1}, \eqref{eq:HA2}, \eqref{eq:HA3} and \eqref{eq:HA1}, \eqref{eq:HA2}, \eqref{eq:HA3} in \eqref{est::finalH} give the result.
\qquad \end{proof}

\begin{remark}
The result obtained is however a little bit disappointing: in the case of the linear driver the result can be improved and it can be shown that the error goes to zero even using an exponential law. In the general case, proposition \ref{propGhost} gives us that the variance is finite using an exponential law for $\rho$ only if  $p \le 2$. 
\end{remark}
\begin{remark}
In the case of non constant coefficients, most of the time it is necessary to use an Euler scheme. 
As ready pointing out in \cite{warin2017variations}, the first estimator has an exploding variance  because the integration by part has to be achieved on the first time step of the Euler scheme using a mesh of size $\Delta t$. It gives a Malliavin weight  in $O(\frac{1}{\sqrt{\Delta t}})$ leading to an explosion in variance as the step size goes to zero. This second estimator doesn't suffer from this problem.
\end{remark}

\subsection{Numerical tests for the semi-linear case}
In this section we give some numerical results illustrating the previous results obtained.
In the whole section the number of particles taken at each level will be a sequence $(N_i^{ipart})_{i\ge 0}$ indexed by $ipart$ such that:
\begin{align}
\label{repPart}
N_i^{ipart}= N_i^0 \times 2^{ipart}.
\end{align}

\subsection{A B\"urgers test case.}
\begin{figure}[h!]
\begin{minipage}[b]{0.49\linewidth}
  \centering
 \includegraphics[width=\textwidth]{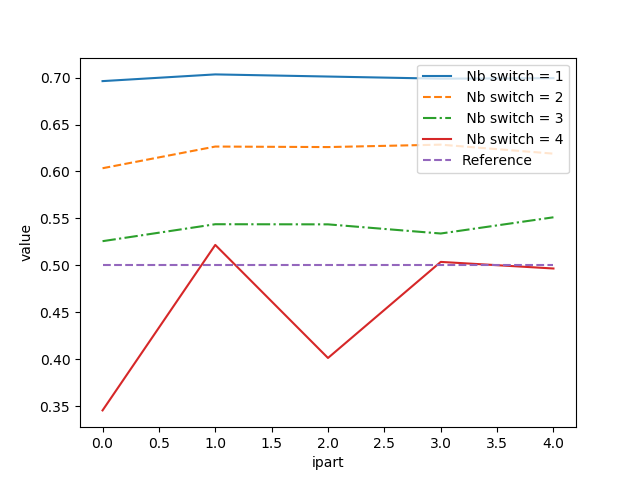}
 \caption*{$\lambda=0.1$.}
 \end{minipage}
\begin{minipage}[b]{0.49\linewidth}
  \centering
 \includegraphics[width=\textwidth]{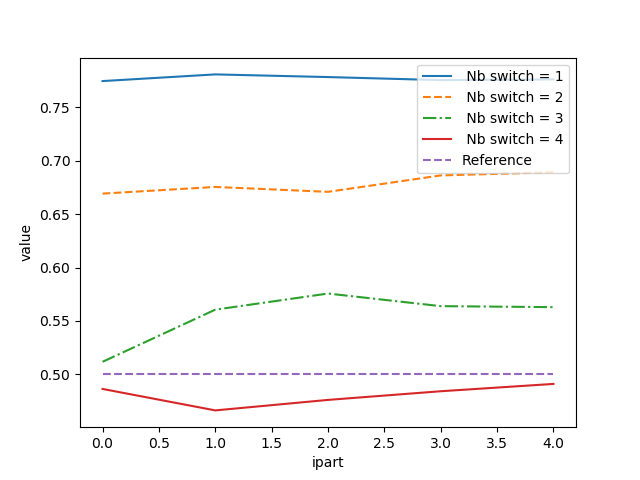}
 \caption*{$\lambda=0.2$.}
 \end{minipage}
 \caption{\label{fig:case3NoGhostGamD10} B\"urgers case: convergence with estimator \eqref{eq:estimSemi1} in dimension 10 using  a gamma law with $u=0.8$.  }
 \end{figure}

\begin{figure}[h!]
\begin{minipage}[b]{0.49\linewidth}
  \centering
 \includegraphics[width=\textwidth]{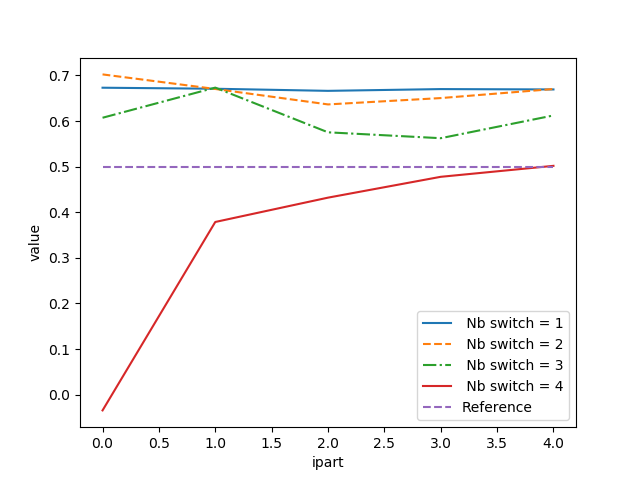}
 \caption*{$\lambda=0.1$}
 \end{minipage}
\begin{minipage}[b]{0.49\linewidth}
  \centering
 \includegraphics[width=\textwidth]{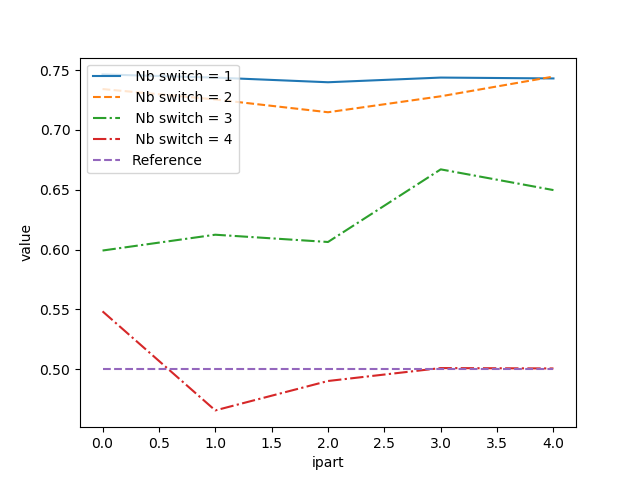}
 \caption*{$\lambda=0.2$}
 \end{minipage}
 \caption{\label{fig:case3NoGhostGamD20} B\"urgers case: convergence with estimator \eqref{eq:estimSemi1} in dimension 20 using  a gamma law with $u=0.8$.  }
 \end{figure}

We take the test case proposed in \cite{weinan2017deep} which is derived from a test in \cite{chassagneux2014linear}.
We take the same parameters as in \cite{weinan2017deep}:
$\mu=0$,
$\sigma = d \I_d$, $T=1$, 
the driver is given for $x \in\R^d, y\in \R, z \in \R^d$ by
\begin{align*}
f(t,x,y,z) = (y- \frac{2+d}{2d}) (d \sum_{i=1}^d z_i),
\end{align*}
and the final function is :
\begin{align*}
g(x) =\frac{e^{T+\frac{1}{d} \sum_{i=1}^d x_i}}{1+e^{T+\frac{1}{d} \sum_{i=1}^d x_i}}.
\end{align*}
The explicit solution given by \cite{weinan2017deep}  is
\begin{align*}
u(t,x) =\frac{e^{t+\frac{1}{d} \sum_{i=1}^d x_i}}{1+e^{t+\frac{1}{d} \sum_{i=1}^d x_i}}.
\end{align*}
We solve the problem in dimension $d=10$ and $d=20$ at date $t=0$ and for $x= 0 \un_d$ such that the reference is equal to $0.5$.\\
We test the different schemes for a gamma law  given by \eqref{rho} with $\lambda=0.1$ and $\lambda=0.2$.\\
For the first scheme we take $u=0.8$ in \eqref{rho}  and give the results obtained on figure \ref{fig:case3NoGhostGamD10} and \ref{fig:case3NoGhostGamD20} for estimator \eqref{eq:estimSemi1}
with  $(N_0^0, N_1^0, N_2^0, N_3^0)  = (1000,40,40,30)$ in  \eqref{repPart}.\\
In dimension 10 and 20, we obtain  very good results with 4 switches and $ipart=4$,
\begin{itemize}
\item getting in dimension $10$ a value $0.496$  for $\lambda=0.1$ in  $80$ seconds and $0.4910$  for $\lambda=0.2$ in  $1000$ seconds,
\item getting  in dimension $20$ a value $0.501$  for $\lambda=0.1$ in  $350$ seconds and $0.5006$  for $\lambda=0.2$ in  $1400$ seconds.
\end{itemize}

\begin{figure}[h!]
\begin{minipage}[b]{0.49\linewidth}
  \centering
 \includegraphics[width=\textwidth]{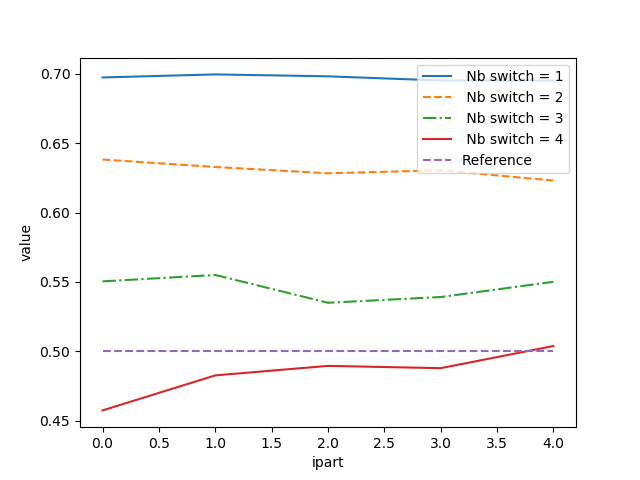}
 \caption*{$\lambda=0.1$.}
 \end{minipage}
\begin{minipage}[b]{0.49\linewidth}
  \centering
 \includegraphics[width=\textwidth]{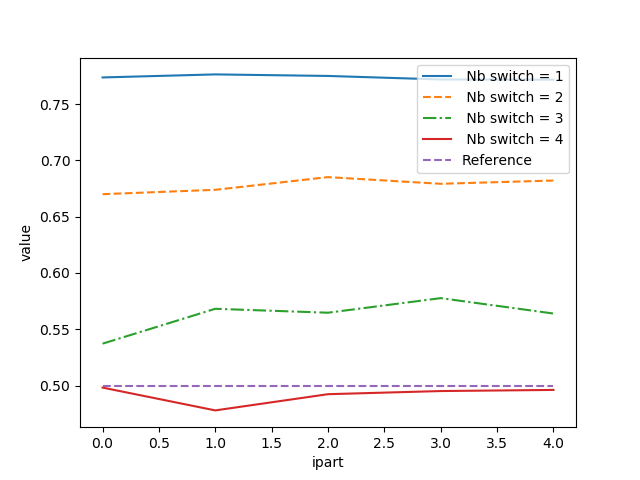}
 \caption*{$\lambda=0.2$.}
 \end{minipage}
 \caption{\label{fig:case3GhostGamD10} B\"urgers case: convergence with estimator \eqref{eq:estimSemi2} in dimension 10 using a gamma law with $u=0.9$.  }
 \end{figure}

\begin{figure}[h!]
\begin{minipage}[b]{0.49\linewidth}
  \centering
 \includegraphics[width=\textwidth]{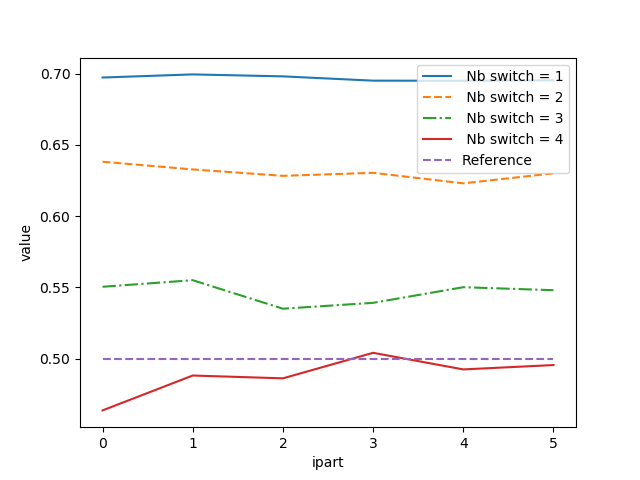}
 \caption*{$\lambda=0.1$.}
 \end{minipage}
\begin{minipage}[b]{0.49\linewidth}
  \centering
 \includegraphics[width=\textwidth]{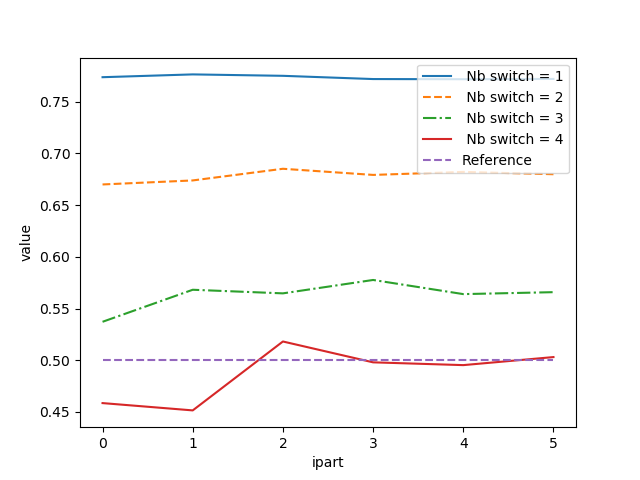}
 \caption*{$\lambda=0.2$.}
 \end{minipage}
 \caption{\label{fig:case3GhostGamD20} B\"urgers case: convergence with estimator \eqref{eq:estimSemi2} in dimension 20 with a gamma law with $u=0.9$.  }
 \end{figure}

\begin{figure}[h!]
\begin{minipage}[b]{0.49\linewidth}
  \centering
 \includegraphics[width=\textwidth]{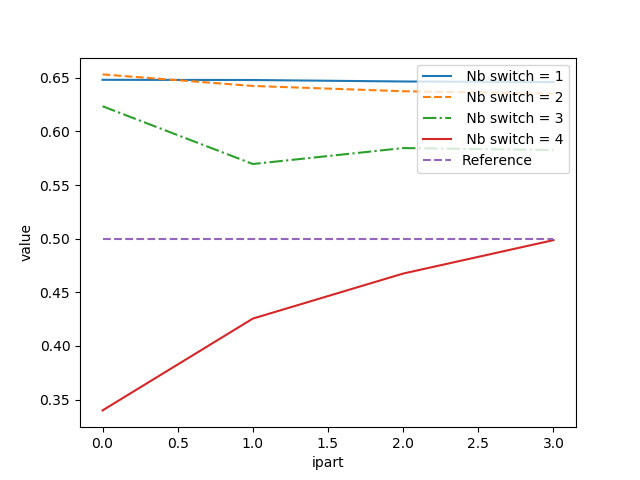}
 \caption*{$\lambda=0.1$.}
 \end{minipage}
\begin{minipage}[b]{0.49\linewidth}
  \centering
 \includegraphics[width=\textwidth]{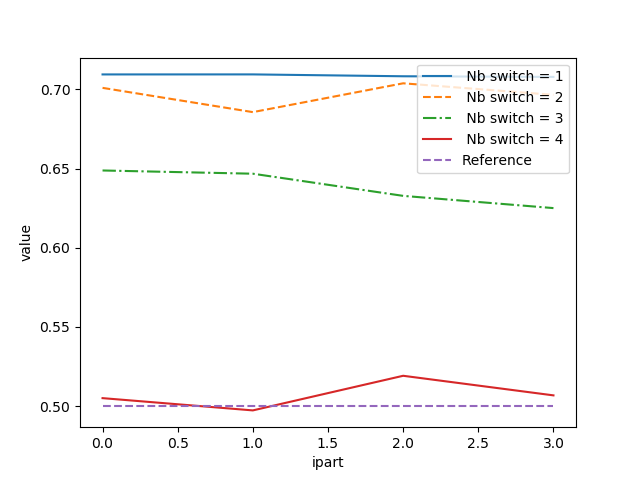}
 \caption*{$\lambda=0.2$.}
 \end{minipage}
 \caption{\label{fig:case3GhostExpD20} B\"urgers case: convergence with estimator \eqref{eq:estimSemi2} in dimension 20 with an exponential law.  }
 \end{figure}

For the second scheme \eqref{eq:estimSemi2}, taking $u=0.9$ in 
\eqref{rho}, we give  the results obtained on figure \ref{fig:case3GhostGamD10} and \ref{fig:case3GhostGamD20} for estimator \ref{eq:estimSemi1}
with  $(N_0^0, N_1^0, N_2^0, N_3^0)  = (1000,40,40,30)$ in  \eqref{repPart}.
In dimension 10 and 20 , we obtain  very good results with 4 switches and $ipart=3$,
\begin{itemize}
\item getting in dimension $10$ a value $0.487$  for $\lambda=0.1$ in  $27$ seconds and $0.49509$  for $\lambda=0.2$ in  $130$ seconds,
\item getting  in dimension $20$ a value $0.5040$  for $\lambda=0.1$ in  $40$ seconds and $0.4979$  for $\lambda=0.2$ in  $186$ seconds.
\end{itemize}
 At last we use an exponential law for $\rho$, leading to $u=1$ in \eqref{rho} and we use the  estimator \eqref{eq:estimSemi2}. We take  $(N_0^0, N_1^0, N_2^0, N_3^0)  = (1000,40,40,4)$ in  \eqref{repPart} and only give the results on figure \ref{fig:case3GhostExpD20} for the most difficult case $d=20$.\\
 With $4$ switches, $ipart=3$ we get very good  results :
 \begin{itemize}
 \item 
 A solution equal to $0.499$  with a computational time equal to $14$ seconds with $\lambda=0.1$,
 \item A solution equal to $0.506$  with a computational time equal to  $55$ seconds with $\lambda=0.2$. 
 \end{itemize}
 The variance using exponential laws seems to be lower and the generation of an exponential law takes far less time than with general gamma laws.
 
 \subsubsection{A second case}
We then take the HJB equation test case taken from \cite{weinan2017deep},\cite{han2017overcoming}, \cite{chassagneux2016numerical}.
As in \cite{han2017overcoming}  we solve the problem in dimension 100 to show the efficiency with the same  characteristics as in \cite{weinan2017deep},\cite{han2017overcoming}: 
\begin{align*}
  \mu= & 0, \quad \sigma = \sqrt{2} \I_d,  \quad  T=1, \\
f(t,x,z) =& - \theta ||z||^2_2,
\end{align*}
such that a semi-explicit solution is available :
\begin{align}
\label{eq:semiAnal}
u(t,x) =- \frac{1}{\theta} \log \big(  \E[ e^{ - \theta g(x + \sqrt{2} W_{T-t})}]\big).
\end{align}
In the example, we take the $g$ function as in \cite{weinan2017deep}:
\begin{align*}
g(x)= \log(\frac{1 + ||x||_2^2}{2} ),
\end{align*}
and we want to estimate the solution at date $t=0$ and for $x=  0\un_d$  
using our algorithm. \\
We treat three cases with increasing  difficulty by taking $\theta =1$, then $\theta=10$ then at last $\theta=20$. The difficulty comes from an increasing value of the non linearity.\
Using a Monte Carlo approximation of equation \eqref{eq:semiAnal} we get some references and a good approximation of the solution
is $4.59$ with $\theta=1$,  $4.49$ with $\theta=10$, and $4.36$ with $\theta=20$.
In order to fit the framework we modify the non linearity to 
\begin{align*}
f(t,x,z) = - \theta  \min(||z||^2_2, 1),
\end{align*}
but the truncation has in fact no effect on the method.\\

\begin{figure}[h!]
\begin{minipage}[b]{0.49\linewidth}
  \centering
 \includegraphics[width=\textwidth]{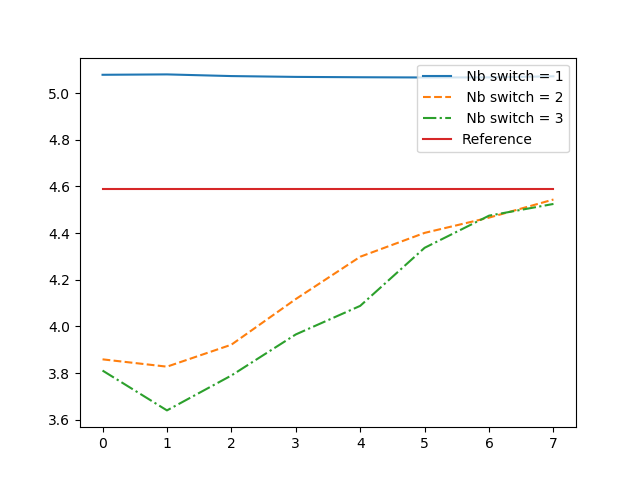}
 \caption*{$\lambda=0.1$}
 \end{minipage}
\begin{minipage}[b]{0.49\linewidth}
  \centering
 \includegraphics[width=\textwidth]{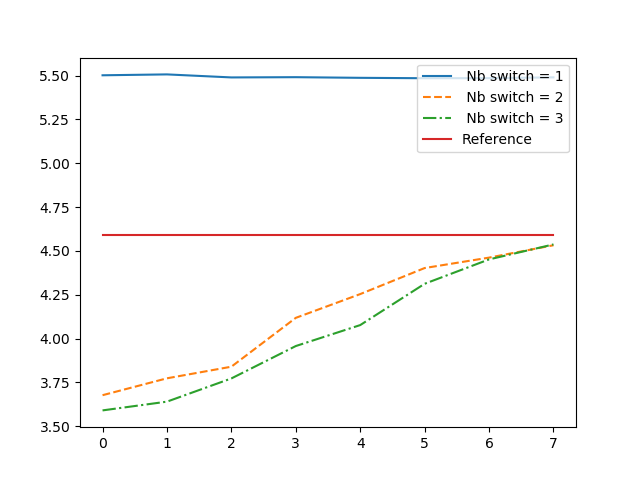}
 \caption*{$\lambda=0.2$}
 \end{minipage}
 \caption{\label{fig:case3NoGhostGamTheta1} HJB convergence case for $\theta=1$ and estimator \ref{eq:estimSemi1} using a gamma law with $u=0.8$.}
\end{figure}

First for $\theta=1$ we plot  on figure \ref{fig:case3NoGhostGamTheta1} the results obtained by estimator \eqref{eq:estimSemi1} taking in  equation \eqref{repPart}:
$(N_0^0, N_1^0, N_2^0)  = (1000,20,20)$.
The convergence with this scheme is quite slow for both $\lambda=0.1$ and $\lambda=0.2$. The number of switches to take to have a very good accuracy seems to be equal to 3 but
the variance observed with this estimator is quite high. Experiments with $\theta=10$ or $\theta=20$ show that this kind of singularity is hard to cope with so high Lipschitz constants of the driver.
\\

\begin{figure}[h!]
\begin{minipage}[b]{0.49\linewidth}
  \centering
 \includegraphics[width=\textwidth]{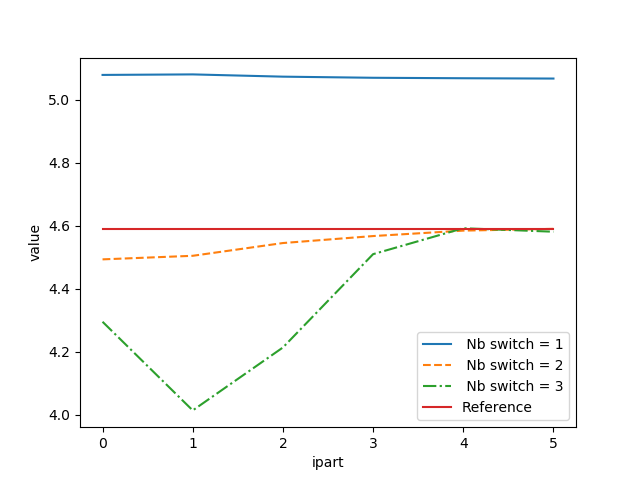}
 \caption*{$\lambda=0.1$}
 \end{minipage}
\begin{minipage}[b]{0.49\linewidth}
  \centering
 \includegraphics[width=\textwidth]{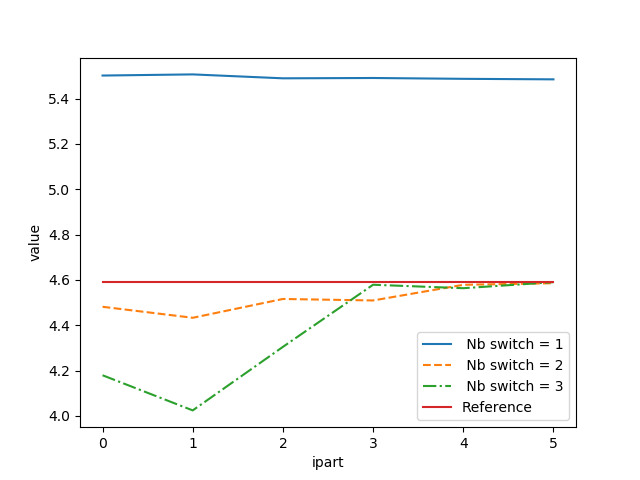}
 \caption*{$\lambda=0.2$}
 \end{minipage}
 \caption{\label{fig:case3GhostGamTheta1} HJB convergence case for $\theta=1$ and estimator \eqref{eq:estimSemi2} using a gamma law with $u=0.9$.}
\end{figure}

\begin{figure}[h!]
\begin{minipage}[b]{0.49\linewidth}
  \centering
 \includegraphics[width=\textwidth]{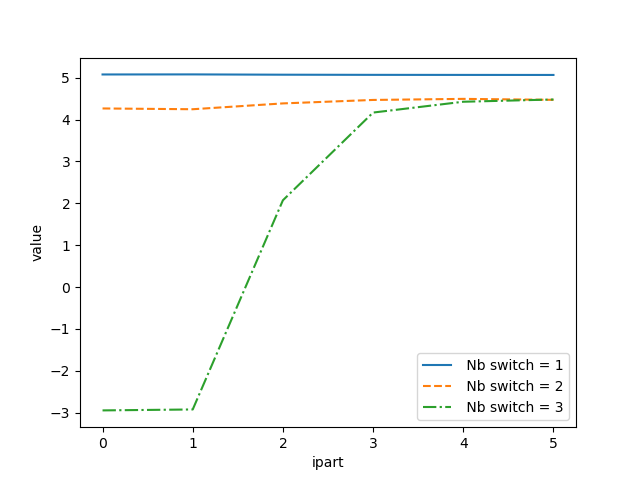}
 \caption*{$\lambda=0.1$}
 \end{minipage}
\begin{minipage}[b]{0.49\linewidth}
  \centering
 \includegraphics[width=\textwidth]{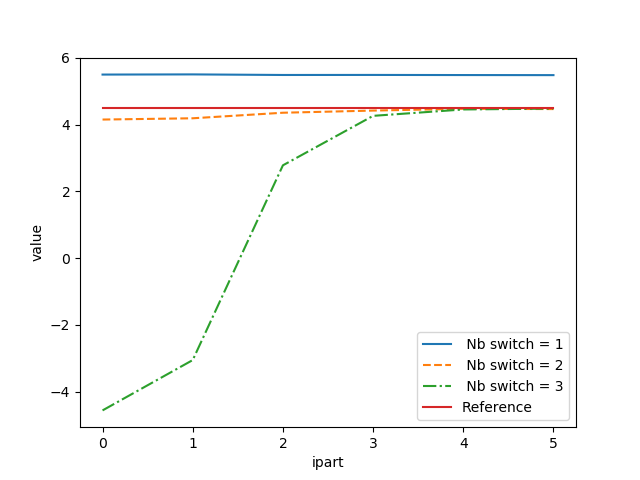}
 \caption*{$\lambda=0.2$}
 \end{minipage}
 \caption{\label{fig:case3GhostGamTheta10} HJB convergence case for $\theta=10$ and estimator \eqref{eq:estimSemi2} using  a gamma law with $u=0.9$.}
\end{figure}

\begin{figure}[h!]
\begin{minipage}[b]{0.49\linewidth}
  \centering
 \includegraphics[width=\textwidth]{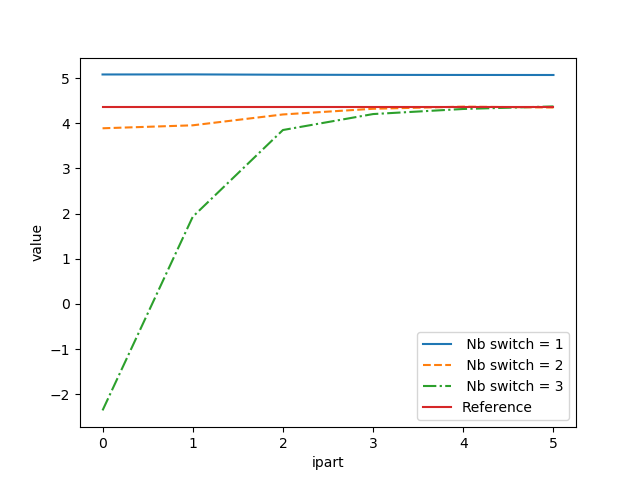}
 \caption*{$\lambda=0.1$}
 \end{minipage}
\begin{minipage}[b]{0.49\linewidth}
  \centering
 \includegraphics[width=\textwidth]{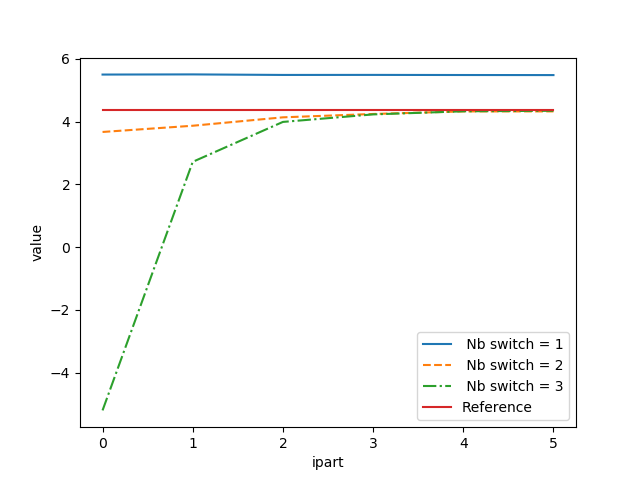}
 \caption*{$\lambda=0.2$}
 \end{minipage}
 \caption{\label{fig:case3GhostGamTheta20} HJB convergence case for $\theta=20$ and estimator \eqref{eq:estimSemi2} using  a gamma law with $u=0.9$.}
\end{figure}

\begin{figure}[h!]
\begin{minipage}[b]{0.49\linewidth}
  \centering
 \includegraphics[width=\textwidth]{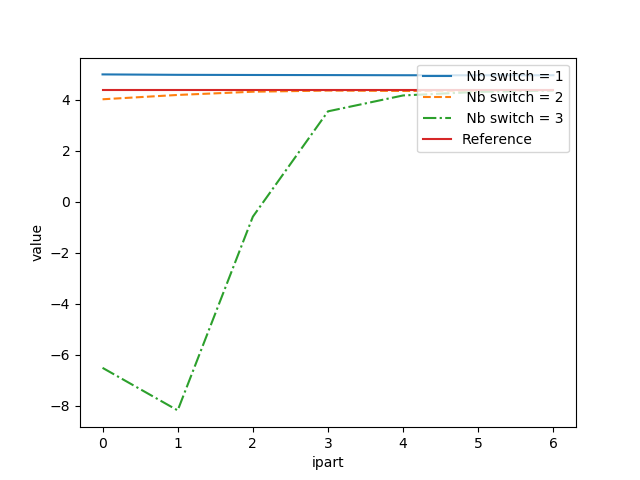}
 \caption*{$\lambda=0.1$}
 \end{minipage}
\begin{minipage}[b]{0.49\linewidth}
  \centering
 \includegraphics[width=\textwidth]{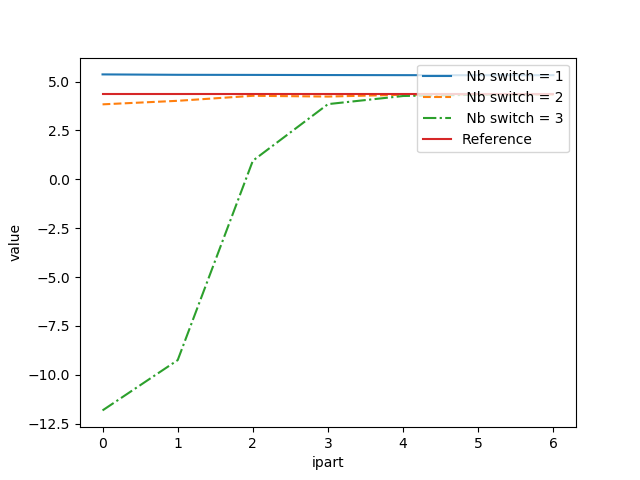}
 \caption*{$\lambda=0.2$}
 \end{minipage}
 \caption{\label{fig:case3GhostExpTheta20} HJB convergence case for $\theta=20$ and estimator \eqref{eq:estimSemi2} using an exponential law.}
\end{figure}

For the three values of $\theta$ on figures
\ref{fig:case3GhostGamTheta1}, \ref{fig:case3GhostGamTheta10} and \ref{fig:case3GhostGamTheta20}
we plot the results for gamma law given by \eqref{rho} with $u=0.9$ and estimator \eqref{eq:estimSemi2}  taking in  equation \eqref{repPart}:
\begin{itemize}
\item $(N_0^0, N_1^0, N_2^0)  = (1000,10,1)$ for $\theta=1$,
\item $(N_0^0, N_1^0, N_2^0)  = (1000,10,5)$ for $\theta=10$,
\item $(N_0^0, N_1^0, N_2^0)  = (1000,40,20)$ for $\theta=20$,
\end{itemize}
such that as the different Lipschitz  constants increase we increase the numbers of samples taken in inner nesting.\\
As for the results obtained using estimator \eqref{eq:estimSemi2},  a good accuracy is obtained using 2 
switches:
\begin{itemize}
\item with $\theta=1$, $\lambda=0.1$, $ipart =3$, we get $4.57$ whereas taking $\lambda=0.2$, $ipart =4$ gives $4.58$.
\item with $\theta=10$,  $\lambda=0.1$, $ipart =3$, we get $4.47$ whereas taking $\lambda=0.2$, $ipart =4$ gives $4.48$.
\item with $\theta =20$, $\lambda=0.1$, $ipart =4$, we get $4.37$ whereas taking $\lambda=0.2$, $ipart =4$ gives $4.33$.
\end{itemize}
At last on figure \ref{fig:case3GhostExpTheta20}, we plot the solution  obtained  for the most difficult case  ($\theta=20$)  using an exponential law for $\rho$ and $(N_0^0, N_1^0, N_2^0)  = (1000,40,10)$.
A good solution is obtained in $10$ seconds for $\lambda=0.1$ with 2 switches with a precision of less than $0.5\%$ using $(N_0, N_1)= (8000,  320)$.
To get a very accurate solution  with a precision of $0.1\%$ with both $\lambda$  it is then necessary  to use 3 switches with $(N_0,N_1,N_2) = (64000,  2560, 640)$ and the
computational time explodes to nearly $30000$ seconds with $\lambda=0.1$ and $80000$ seconds with $\lambda=0.2$.

\section{Conclusion}
An effective method to solve semi-linear equations has been developed and tested. The most effective way to solve these equations consists in taking the second scheme proposed to treat the gradient term.
The limitation due to Lipschitz constant and the maturity of the problem can be easily postponed using cluster of CPU or perhaps GPU.\\
Besides even if cannot prove that the second scheme can be used with an exponential law, it seems to be the most effective. 
A better understanding of its efficiency could pave the way to solve full non linear PDEs.